\documentstyle[amsfonts,11pt]{article}

\newtheorem{theorem}{Theorem}[section]
\newtheorem{lemma}{Lemma}[section]
\newtheorem{proposition}{Proposition}[section]
\newtheorem{remark}{Remark}[section]

\begin{document}

\author{Ion Grama$^1$ and Michael Nussbaum \\
{\it Institute of Mathematics, Chi\c{s}in\u{a}u, Moldova and }\\
{\it Weierstrass Institute, Berlin, Germany}}
\title{{\bf Asymptotic Equivalence for Nonparametric}
       {\bf Generalized Linear Models }}
\date{November 1996\\
Revised: January 1998}
\maketitle

\begin{abstract}
We establish that a non-Gaussian nonparametric regression model is
asymptotically equivalent to a regression model with Gaussian noise. The
approximation is in the sense of Le Cam's deficiency distance $\Delta $; the
models are then asymptotically equivalent for all purposes of statistical
decision with bounded loss. Our result concerns a sequence of independent
but not identically distributed observations with each distribution in the
same real-indexed exponential family. The canonical parameter is a value $%
f(t_i)$ of a regression function $f$ at a grid point $t_i$ (nonparametric
GLM). When $f$ is in a H\"{o}lder ball with exponent $\beta >\frac 12 ,$ we
establish global asymptotic equivalence to observations of a signal $\Gamma
(f(t))$ in Gaussian white noise, where $\Gamma $ is related to a variance
stabilizing transformation in the exponential family. The result is a
regression analog of the recently established Gaussian approximation for the
i.i.d.\ model. The proof is based on a functional version of the Hungarian
construction for the partial sum process.
\end{abstract}

\begin{table}[b]
\rule{3cm}{0.01cm}\newline
\par
{\small $^1$ Research supported by the Deutsche Akademische Austauschdienst
and by the Deutsche Forschungsgemeinschaft. \newline
 }
\end{table}

\section{Introduction\label{sec:Intro}}

\setcounter{equation}{0}

The remarkable success of the Le Cam's asymptotic theory is mostly due to
the power of the concept of {\it weak convergence} of statistical
experiments, which can be established via LAN conditions. Weak convergence
takes place for experiments localized at the normalizer rate for the
underlying central limit theorem, i.~e.~the usual $n^{-1/2}.$ However this
is useless  if estimators have a slower rate of convergence, which happens with
 nonparametric estimation problems of  the ''ill posed'' type.
It is therefore  natural to abandon the localization concept in this case
and to replace limits of experiments by approximations in the sense of Le
Cam's deficiency pseudodistance $\Delta .$ The $\Delta $-distance can be
accessed via coupling of likelihood processes and new results on strong
approximation for sums of random variables. We refer to Koltchinskii \cite
{Kol} for a result on empirical processes in the i.~i.~d.~case.

The global $\Delta $-distance for nonparametric Gaussian experiments was
first studied by Brown and Low \cite{Br-Low}, who showed that a normal
nonparametric regression is asymptotically equivalent to its continuous
version - the signal recovery problem in Gaussian white noise. Then in
Nussbaum \cite{Nuss} it was established that density estimation from i.~i.~d.~data
 on an interval is asymptotically equivalent to the signal recovery
problem, where the signal is the root density. The two sequences of
experiments are then {\it accompanying }in the sense that their deficiency
pseudodistance $\Delta $ tends to zero. This can be regarded as the natural
generalization of the classical local asymptotic normality theory to "ill
posed" problems. The implication for decision theory is "automatic" transfer
of risk bounds from one sequence to another.

The purpose of the present paper is to accomplish a logical next step in
these developments, i.~e.~to treat the case of {\it non-Gaussian
nonparametric regression. } Our model is such that at points $t_i=i/n,$ $%
i=1,...,n,$ we observe independent r.v.'s $X_i,$ which follow a distribution
from an exponential family ${\cal P}$ with parameters $\theta _i=f(t_i)\in
\Theta ,$ where $f:[0,1]\rightarrow \Theta $ is an unknown function to
be estimated. The function $f$ is assumed to belong to the smoothness class $%
\Sigma .$ The main result of the paper is asymptotic equivalence of this
model to a Gaussian experiment of the {\it homoscedastic} form
\begin{equation}
dY_t^n=\Gamma \left( f(t)\right) dt+\frac 1{\sqrt{n}}dW_t,\quad t\in [0,1],
\label{I.1}
\end{equation}
with $f\in \Sigma ,$ where $\Gamma (\theta ):\Theta \rightarrow R$ is a
function such that $\Gamma ^{\prime }(\theta )=I(\theta )^{1/2}$ and $%
I(\theta )$ is the Fisher information in the local exponential family ${\cal %
P}$ (see Sections \ref{sec:VST} and \ref{sec:VS} for a relation with the
variance-stabilizing transformation pertaining to ${\cal P}$). Here $W$ is the
standard Wiener process, and $f$ runs a set of functions in a H\"{o}lder
ball with exponent $\beta >\frac 12$.

Note that our function $f$ is tied to the canonical parametrization of the
exponential family, while the ''natural'' parameter (the intensity for the
Poisson case etc.) is generally different. But there is a smooth parameter
transformation $\lambda =b(\theta )$ (defined in Section \ref{sec:VST}
below) which permits to formulate global results of the type (\ref{I.1}) in
''natural'' regression models. Some examples are:

\begin{description}

\item[{\bf [1] Poisson case:} ]   $X_i$ is Poisson$[g(t_i)],$ where $g$ is a function
on $[0,1]$ in a H\"{o}lder ball with exponent $\beta >\frac 12,$ with values
in $[\epsilon ,\epsilon ^{-1}],$ for some $\epsilon >0$. The Poisson
intensity is $\lambda =b(\theta ),$ where $b$ is a strictly increasing
smooth function. Defining the function $F(\lambda )=\Gamma (b^{-1}(\lambda
)),$ we obtain (see Section \ref{sec:Examples} below) $F(\lambda )=2\sqrt{%
\lambda }$. The accompanying Gaussian experiment is
\begin{equation}
dY_t^n=2\sqrt{g(t)}\,dt+\frac 1{\sqrt{n}}dW_t,\quad t\in [0,1].  \label{I.e2}
\end{equation}

\item[{\bf [2] Bernoulli case:}]   $X_i$ is Binomial$[1,g(t_i)],$ $g$ as above, but
with values in $[\epsilon ,1-\epsilon ],$ for some $\epsilon >0$. The
natural parameter is $\lambda =b(\theta ),$ for some function $b$ with
properties as above. We have $F(\lambda )=2\arcsin (\sqrt{\lambda }),$ and
the accompanying Gaussian experiment is
\[
dY_t^n=2\arcsin \sqrt{g(t)}\,dt+\frac 1{\sqrt{n}}dW_t,\quad t\in [0,1].
\]

\item[{\bf [3] Gaussian variance case:} ]   $X_i$ is $N[0,g(t_i)],$ $g$ as in example
1. We have $F(\lambda )=2^{-1/2}\log \lambda ,$ and
\begin{equation}
dY_t^n=\frac 1{\sqrt{2}}\log g(t)\,dt+\frac 1{\sqrt{n}}dW_t,\quad t\in [0,1].
\label{I.e1}
\end{equation}

\end{description}

For more details see Section \ref{sec:Examples}. We chose to give here the
continuous versions of the accompanying experiments, but the discrete
versions are also available.

The motivation for this paper can be concisely expressed in three points.

\begin{description}

\item [ {\bf [A] }]  The first example essentially recovers the result of \cite{Nuss}
on i.i.d.\ data with density $f$ on the unit interval, since the proof in
\cite{Nuss} used a Poisson approximation. In the second example, $F(\lambda
)=2\arcsin (\sqrt{\lambda })$ is the well known variance stabilizing
transformation for the binomial distribution. For scale parameters like the
Gaussian variance, the logarithm is also known to be a stabilizing
transformation, and the same applies to $2\lambda ^{1/2}$ in the Poisson
case. Thus, in the evolving theory of asymptotic equivalence of experiments,
we have achieved a better understanding of where global closed form
approximations like (\ref{I.1}) arise from. The formal connection to
differential geometric theory in statistics seems very interesting and
remains to be explored (cp.\v {C}encov \cite{Cen}, Amari et al. \cite{Amar}).

\item [ {\bf [B] }]  The case when the function $f$ is in a linear parametric class $%
\left\{ f(x)=\beta x,\;\beta \in R\right\} $ is known as a {\em generalized
linear model} (GLM). The inverse of the transformation $b(\theta )$ would
then be the canonical link function. Accordingly, our model is of the type
{\em nonparametric GLM}, cp. Green and Silverman \cite{Green-Silv}, sec.
5.1.2. Models like these, which offer tremendous flexibility, have received
much attention in the recent literature, see also Fan and Gijbels \cite
{Fan-Gijb}. It would be beyond the scope of this paper to treat the many
(semiparametric) variants and extensions; we refer to \cite{Green-Silv} and
\cite{Fan-Gijb}. In particular we do not discuss logit and probit analysis
in our context (cp. example 3, p. 92 in \cite{Green-Silv}). Empirical
process theory has also been applied by Mammen and van de Geer \cite
{Mamm-Geer} in our model (in a more general variant), for constructing
estimators.

\item [ {\bf [C] }]  There are implications for time series models. Example 3 leads on
to the white noise equivalence for the {\it spectral density model} for a
Gaussian stationary sequence; cf. Golubev and Nussbaum \cite{Gol}.
Furthermore, Example 3 is related to discrete observations of a diffusion
process:
\[
dY_t^n=g^{1/2}(Y_t^n,t)dW_t,\quad Y_0^n=y_0,\quad t\in [0,1].
\]
Suppose observations in points $t_i=i/n,$ $i=1,\ldots ,n,$ where $g$ is
unknown. Nonparametric estimation of $g$ has recently been considered by
Genon-Catalot, Laredo and Picard \cite{GCLP}, Florens-Zmirou \cite{FZ}.
Example 3 might be seen as a possible pilot result for those models, where
the distributions of processes on $[0,1]$ are mutually orthogonal and the
asymptotics is given by  grid refinement.

\end{description}

The standard method of proof is to establish first a local version of (\ref
{I.1}) and then to globalize it by means of a preliminary estimator. We
obtain our initial local approximation in the {\it heteroscedastic} form
\begin{equation}
dY_t^n=f(t)\,dt+\frac 1{\sqrt{n}}I\left( f_0(t)\right) ^{-1/2}dW_t,\quad
t\in [0,1],  \label{I.2}
\end{equation}
where $f$ is in a shrinking neighborhood of a function $f_0$ and $I(\theta )$
is the Fisher information in the local exponential family ${\cal P}$ (given
in its canonical form). To obtain a global asymptotic equivalence, the
function $f_0,$ which was technically assumed ''known'', has to be replaced
by a preliminary estimator. However, the homoscedastic form (\ref{I.1}) can
be obtained only if the function $f_0$ does not show up explicitly in the
local approximation. The problem arises to find a transformation $\Gamma
(\theta )$ on the parameter space of the local exponential family ${\cal P}$
such that an asymptotically equivalent form of (\ref{I.2}) would be (\ref
{I.1}), with $f$ in a neighborhood of $f_0.$ It is easy to show that such a
transformation exists, for any exponential family ${\cal P}.$ This question
is related to that of a variance-stabilizing transformation. Indeed, for the
Poisson observations of example [1], one can prove easily that an
accompanying local experiment, besides (\ref{I.2}), is also given by
\begin{equation}
dY_t^n=g(t)\,dt+\frac 1{\sqrt{n}}\sqrt{g_0(t)}dW_t,\quad t\in [0,1],
\label{I.3}
\end{equation}
with $g$ in a neighborhood of $g_0$ (for an analogy with the i.~i.~d.~model
see also Nussbaum \cite{Nuss}). The observations in (\ref{I.3}) have roughly
a Poisson character - the expectation is $g(t)$ and the variance is
approximately $g(t)$ as well (since $g$ is in a neighborhood of $g_0$). For the
Poisson distribution, the square root is a variance stabilizing
transformation, which agrees with (\ref{I.e2}).

For the proof of the local heteroscedastic approximation (\ref{I.2}) we
establish a functional Hungarian construction for the partial sums of
independent but nonidentically distributed random variables. The result is
similar to that of Koltchinskii \cite{Kol} for the empirical process, but
the generality of assumptions ( non-identity and non-smoothness of distributions
 of the summands) substantially complicates the problem. This is treated separately in
\cite{Gr-Nu} based on methods developed in Sakhanenko \cite{Sakh}. Due to
the particularly simple structure of our nonparametric exponential model, we
can straightforwardly apply our strong approximation result to couple the
likelihood process with that of an appropriately chosen Gaussian experiment.
Again {\em coupling of likelihood processes} is the key idea for proving
asymptotic equivalence, as in \cite{Nuss}. However we would like to mention
that our KMT result is useful in more general situations also.

An essential step in proving a local approximation result like (\ref{I.2})
is to study the local experiments generated by the fragments of observations
$X_i$ over shrinking time intervals. These experiments we call {\it doubly
local}. With an appropriate choice of the length of the shrinking interval
and after rescaling it to the unit interval, a doubly local experiment can
be viewed as a local experiment on the interval $[0,1]$ of the usual type, but
now with a neighborhood of the ''almost root-$n$'' size $(n/\log n)^{-1/2}.$
A similar renormalization technique is known to be effective for pointwise
estimation in nonparametrics, cf. Donoho and Liu \cite{Don-Liu} and Low \cite
{Low}. We also refer to Millar \cite{Mill} for $n^{-1/2}$ -shrinking
neighborhoods in the context of nonparametric estimation. Were it not for
the $\log $-factor in $(n/\log n)^{-1/2},$ these rescaled experiments on the
interval $[0,1]$ would converge to a Gaussian limit in the sense of $\Delta $%
. The motivation for applying the Hungarian construction at this stage is,
roughly speaking, to obtain a good rate for this convergence, i.~e.~$%
n^{-1/2} $ up to some logarithmic factor. We thus implicitly address a
question of Le Cam on rates for convergence of experiments (cp.~the remark
on p.~509 of \cite{LeCam}). The use of the Hungarian construction for this
purpose is in line with its original motivation, i.~e.~optimal rates in the
functional central limit theorem for the partial sum process.

Our results are formulated in terms of $\Delta $-distance approximations; we
do not exhibit the recipes (Markov kernels) which transfer a decision
function in one experiment to the other. The problem of constructive
equivalence is an important issue and some promising research in this
direction is going on (Brown and Low \cite{Br-Low-BER}). Markov kernels in
the present case can be extracted from the Hungarian construction, but this
is beyond of the scope of the present paper.

The method of the proof is similar to that in Nussbaum \cite{Nuss}. We
utilize the natural independence structure of our regression model for
decomposing it into fragments on appropriate shrinking intervals (or
equivalently into doubly local experiments). We then apply the functional
version of the KMT construction for establishing a Gaussian local
approximation for the fragments. It is important is to have this
approximation with a rate of convergence that is enough to ''beat'' the
number of shrinking intervals  into which the whole interval $[0,1]$ has
been split. Having obtained the above rate, we still have to ensure the
global approximation on $[0,1]$. This we achieve by passing from the Le Cam
pseudo-distance between statistical experiments to the Hellinger distance
between the corresponding probability measures. The passage is made possible
by a construction of the likelihood processes of the local experiments on
the same probability space with a Gaussian likelihood process. The reason to
work with the Hellinger distance rather than with total variation distance
is the convenient behavior of the former under the product operation for  probability
measures [see (\ref{SE.2})]. This allows to patch together the doubly local
experiments for obtaining the Gaussian approximation (\ref{I.2}) valid
globally on the interval of observations $[0,1],$ but still around a
specified regression function $f_0.$

After that, we choose the variance-stable version of the Gaussian local
approximation as a starting point for globalization over the parameter space
$\Sigma .$ The result is the global approximation (\ref{I.1}). We trace the
rates of convergence throughout, so that the rate of the deficiency distance
approximation can be made explicit.


\section{Background\label{sec:BKG}}

\setcounter{equation}{0}


\subsection{Exponential families of distributions\label{sec:EXP-FAM}}

We will consider a one-dimensional linearly indexed exponential family (see
Brown \cite{Brown} or Le Cam \cite{LeCam}, p. 144), which is described by
means of the following objects:

\begin{itemize}
\item  A measurable space $(X,{\cal B}(X),\mu )$ equipped with the positive
measure $\mu (dx),$ where $X$ is a Borel measurable subset of the real line $%
R$ and ${\cal B}(X)$ is the Borel $\sigma $-field on $X;$

\item  a measurable map $U(x):X\rightarrow R;$

\item  an open (possibly infinite) interval $\Theta $ in $R,$ where the
Laplace transform
\[
L(\theta )=\int_X\exp \{\theta U(x)\}\mu (dx)
\]
is finite.
\end{itemize}

Set $V(\theta )=\log L(\theta ).$ Denote by ${\cal P}=\{P_\theta :\theta \in
\Theta \}$ the set of probability measures $P_\theta $ on the space $(X,%
{\cal B}(X),\mu )$ of the form
\begin{equation}
P_\theta (dx)=\exp \{\theta U(x)-V(\theta )\}\mu (dx).  \label{BKG.1}
\end{equation}
We call ${\cal P}$ an {\it exponential family} on the space $(X,{\cal B}%
(X),\mu )$ with parameter set $\Theta .$ The family of measures ${\cal P}$
defines the {\it exponential experiment} ${\cal E}=(X,{\cal B}(X),{\cal P}),$
which will be the background object for constructing our nonparametric
model. In the case when the measures $P_\theta $ are defined by (\ref{BKG.1}%
), we also say that the exponential family ${\cal P}$ (or the exponential
experiment ${\cal E}$) is given in its {\it canonical form}.

It follows from the definition of an exponential family that the function $%
V(\theta )$ is analytic on $\Theta .$ Also note that, for a fixed $\theta
\in \Theta ,$ the function $V(\theta +t)-V(\theta )$ coincides with the
cumulant generating function
\[
G_\theta (t)=\log \int_{X}\exp \{tU(x)\}P_\theta (dx)
\]
of the r.v.\ $Y(\theta )=U(X(\theta )),$ where the r.v.\ $X(\theta )$ has
the distribution in the exponential family ${\cal P}$ with parameter $\theta
.$ This implies that $\frac{d^m}{d\theta ^m}V(\theta )$ is the cumulant of
order $m$ of the r.v.\ $Y(\theta ).$ In particular
\begin{eqnarray}
V^{\prime }(\theta ) &=&E_\theta Y(\theta ),  \label{BKG.3a} \\
V^{\prime \prime }(\theta ) &=&E_\theta Y(\theta )^2-(E_\theta Y(\theta ))^2,
\label{BKG.3b} \\
V^{\prime \prime \prime }(\theta ) &=&E_\theta Y(\theta )^3-3E_\theta
Y(\theta ))^2E_\theta Y(\theta ))+2(E_\theta Y(\theta ))^3.  \label{BKG.3c}
\end{eqnarray}
For this reason we will call the function $V(\theta )$ the {\it cumulant
generating function} associated with the exponential family ${\cal P}$ or
with the exponential experiment ${\cal E}.$

Consider an exponential experiment ${\cal E}.$ The minimal sufficient
statistic in this experiment is the function $U(x),$ $x\in {X.}$ It is easy
to see that the corresponding Fisher information $I(\theta )$ is
\begin{equation}
I(\theta )=V^{\prime \prime }(\theta )=\int_X\frac{(p_\theta ^{\prime }(x))^2%
}{p_\theta (x)}\mu (dx),  \label{BKG.4}
\end{equation}
where $p_\theta (x)=P_\theta (dx)/\mu (dx)=\exp \{\theta U(x)-V(\theta )\}.$
For any $\theta _0\in \Theta $ and $\varepsilon _0>0,$ denote
\[
B(\theta _0,\varepsilon _0)=\{\theta \in \Theta :|\theta -\theta _0|\leq
\varepsilon _0\}.
\]
Throughout the paper we assume that the following conditions hold true.

\begin{itemize}
\item  The Fisher information is positive on $\Theta ,$ i.e.\
\begin{equation}
I(\theta )>0,\quad \theta \in \Theta .  \label{BKG.Ass 0}
\end{equation}

\item  There exists a (possibly infinite) interval $\Theta _0$ in the
parameter set $\Theta $ such that
\begin{equation}
I_{\min }\leq \inf_{\theta \in \Theta _0}I(\theta ),\quad \sup_{\theta _0\in
\Theta _0}\sup_{\theta \in B(\theta _0,\varepsilon _0)}I(\theta )\leq
I_{\max },  \label{BKG.Ass 1}
\end{equation}
where $I_{\max },$ $I_{\min }$ and $\varepsilon _0$ are positive constants
depending only on the family ${\cal P}.$
\end{itemize}

We will see that condition (\ref{BKG.Ass 1}) can easily be verified for all
examples in Section \ref{sec:Examples}.

Let us denote by $\overline{Y}(\theta )$ the sufficient statistic $Y(\theta
)=U(X(\theta ))$ centered under the measure $P_\theta :\overline{Y}(\theta
)=Y(\theta )-E_\theta Y(\theta ).$ The following assertions are almost
trivial.

\begin{proposition}
\label{P BKG1}Assume that condition (\ref{BKG.Ass 1}) holds true. Then, for
any $|t|\leq \varepsilon _0,$
\[
\sup_{\theta \in \Theta _0}E_\theta \exp \{t\overline{Y}(\theta )\}\leq \exp
\{t^2I_{\max }/2\}.
\]
\end{proposition}

{\noindent{\bf Proof.}\ \ } For the proof it is enough to remark that for $%
|t|\leq \varepsilon _0$
\begin{eqnarray*}
E_\theta \exp \{t\overline{Y}(\theta )\} &=&\exp \{V(\theta +t)-V(\theta
)-tV^{\prime }(\theta )\} \\
&=&\exp \{\frac 12t^2I(\theta +\lambda t)\} \\
&\leq &\exp \{\frac 12t^2I_{\max }\},
\end{eqnarray*}
where $0\leq \lambda \leq 1.$ {\mbox{\ \rule{.1in}{.1in}}}

\begin{proposition}
\label{P BKG2}Assume that condition (\ref{BKG.Ass 1}) holds true. Then
\[
\sup_{\theta _0\in \Theta _0}\sup_{\theta \in B(\theta ,\varepsilon
_0/2)}V^{\prime \prime \prime }(\theta )\leq c,
\]
where $c$ is a constant depending only on $I_{\max }$ and $\varepsilon _0.$
\end{proposition}

{\noindent{\bf Proof.}\ \ } This assertion is an easy consequence of
Proposition \ref{P BKG1} and (\ref{BKG.3a}), (\ref{BKG.3b}), (\ref{BKG.3c}).
{\mbox{\ \rule{.1in}{.1in}}}


\subsection{Variance stabilizing transformation\label{sec:VST}}

Let $X_i,$ $i=1,...,n,$ be a sequence of i.~i.~d.~r.~v.'s each with
distribution in the exponential family ${\cal P}$ for the same parameter $%
\theta \in \Theta .$ Let $V(\theta )$ be the cumulant generating function
associated with the exponential family ${\cal P}.$ Set, for brevity,
$S_n=\frac 1n\sum_{i=1}^nX_i$ and $b(\theta )=V^{\prime }(\theta )=E_\theta
X_1,$ $I(\theta )=V^{\prime \prime }(\theta )=\mathop{\rm Var}_\theta X_1.$
According to the central limit theorem, the sequence $\sqrt{n}\left(
S_n-b(\theta )\right) $ converges weakly, under the measure $P_\theta ,$ to
the normal r.v.\ with zero mean and variance $I(\theta ).$ We are interested
in finding a function $F:R\rightarrow R,$ called {\it the
variance-stabilizing transformation,}{\em \ } which reduces the variance of
the limiting normal law to a constant, i.e.\ such that, under the same
measure $P_\theta $%
\begin{equation}
\sqrt{n}\{F(S_n)-F(b(\theta ))\}\stackrel{d}{\rightarrow }N(0,1),
\label{VST.1}
\end{equation}
for all $\theta \in \Theta .$ Such a transformation exists and is given by
the equation
\begin{equation}
F^{\prime }(b(\theta ))=I(\theta )^{-1/2}.  \label{VST.2}
\end{equation}
The straightforward arguments are similar to those in Barndorff-Nielsen and
Cox \cite{Ba-Nie} (p. 37) or in Andersen et al.\ \cite{Ander} (p. 109).
Indeed, it is easy to see by standard reasoning that two sequences of r.v.'s
$\sqrt{n}\left\{ F\left( S_n\right) -F(b(\theta ))\right\} $ and $\sqrt{n}%
F^{\prime }(b(\theta ))\left\{ S_n-b(\theta )\right\} $ are asymptotically
equivalent in the sense that under the probability $P_\theta $
\begin{equation}
\sqrt{n}\left\{ F\left( S_n\right) -F(b(\theta ))\right\} -\sqrt{n}F^{\prime
}(b(\theta ))\left\{ S_n-b(\theta )\right\} \stackrel{d}{\rightarrow }0
\label{VST.3}
\end{equation}
as $n\rightarrow \infty .$ By the central limit theorem, under the measure $%
P_\theta $
\begin{equation}
\sqrt{n}F^{\prime }(b(\theta ))\left\{ S_n-b(\theta )\right\} \stackrel{d}{%
\rightarrow }N(0,1)  \label{VST.4}
\end{equation}
as $n\rightarrow \infty $ if the function $F$ is chosen such that (\ref
{VST.2}) holds true. From (\ref{VST.3}) and (\ref{VST.4}) we immediately
infer the claim (\ref{VST.1}).

Our proof of the variance-stable form of the asymptotically equivalent
Gaussian experiments follows a similar pattern, see Section \ref{sec:P-VS}.


\subsection{Basic facts on statistical equivalence}

Let $P$ and $Q$ be two probability measures on the measurable space $\left(
\Omega ,{\cal F}\right) .$ The Hellinger distance between the probability
measures $P$ and $Q$ is defined as

\begin{equation}
H^2\left( P,Q\right) =\frac 12{\bf E}_\mu \left( (dP/d\mu )^{1/2}-(dQ/d\mu
)^{1/2}\right) ^2,  \label{SE.1}
\end{equation}
where $P$ and $Q$ are absolutely continuous w.r.t. the probability measure $%
\mu $.

Let $P_1,...,P_n$ and $Q_1,...,Q_n,$ be probability measures on $\left(
\Omega ,{\cal F}\right) .$ Set $P^n=P_1\times ...\times P_n$ and $%
Q^n=Q_1\times ...\times Q_n.$ Then (see Strasser \cite{Str}, lemma 2.19)
\begin{equation}
H^2\left( P^n,Q^n\right) \leq \sum_{i=1}^nH^2\left( P_i,Q_i\right) .
\label{SE.2}
\end{equation}

Let ${\cal E}=\left( \Omega ^1,{\cal F}^1,\left\{ P_\theta :\theta \in
\Theta \right\} \right) $ and ${\cal G}=\left( \Omega ^2,{\cal F}^2,\left\{
Q_\theta :\theta \in \Theta \right\} \right) $ be two statistical
experiments with the same parameter set $\Theta .$ Assume that $(\Omega ^1,%
{\cal F}^1{\cal )}$ and $(\Omega ^2,{\cal F}^2{\cal )}$ are complete
separable  metric (Polish) spaces. The deficiency of the experiment ${\cal E}$
with respect to the experiment ${\cal G}$ is defined as
\[
\delta ({\cal E},{\cal G})=\inf \sup_{\theta \in \Theta }\left\| KP_\theta
-Q_\theta \right\| ,
\]
where the $\inf $ is taken over the set ${\cal M}(\Omega ^1,{\cal F}^2)$ of
all Markov kernels $K$ from $(\Omega ^1,{\cal F}^1{\cal )}$ to $(\Omega ^2,%
{\cal F}^2).$ Le Cam's $\Delta $- distance between the experiments ${\cal E}$
and ${\cal G}$ is defined by
\[
\Delta \left( {\cal E},{\cal G}\right) =\max \left\{ \delta ({\cal E},{\cal G%
}),\delta ({\cal G},{\cal E})\right\} .
\]
Let ${\cal E}^n$ and ${\cal G}^n,$ $n=1,2,...,$ be two sequences of
statistical experiments. We say that ${\cal E}^n$ and ${\cal G}^n$ are {\it %
asymptotically equivalent} if
\[
\Delta \left( {\cal E}^n,{\cal G}^n\right) \rightarrow 0,\quad n\rightarrow
\infty .
\]

We will need a relation between Le Cam and Hellinger distances. Let ${\cal E}%
=(\Omega ^1,{\cal F}^1,\{P_\theta :\theta \in \Theta \})$ and ${\cal G}%
=(\Omega ^2,{\cal F}^2,\{Q_\theta :\theta \in \Theta \})$ be two experiments
with the same parameter set $\Theta .$ Assume that there is some point $%
\theta _0\in \Theta $ such that $P_\theta \ll P_{\theta _0}$ and $Q_\theta
\ll Q_{\theta _0}.$ Suppose that there are versions of the likelihood ratios
$\Lambda ^1(\theta )=P_\theta /dP_{\theta _0}$ and $\Lambda ^2(\theta
)=dQ_\theta /dQ_{\theta _0}$ (as processes indexed by $\theta $) on a common
probability space $(\overline{\Omega },\overline{{\cal F}},{\bf P}).$ Then
the $\Delta $-distance satisfies the inequality
\begin{equation}
\Delta \left( {\cal E},{\cal G}\right) \leq \sqrt{2}\sup_{\theta \in \Theta
}H\left( \Lambda ^1(\theta ),\Lambda ^2(\theta )\right) ,  \label{SE.6}
\end{equation}
where the Hellinger distance between likelihood ratios $\Lambda ^1(\theta )$
and $\Lambda ^2(\theta )$ is defined in analogy to the case of probability
measures:
\begin{equation}
H^2\left( \Lambda ^1(\theta ),\Lambda ^2(\theta )\right) =\frac 12{\bf E}%
\left( \sqrt{\Lambda ^1(\theta )}-\sqrt{\Lambda ^2(\theta )}\right) ^2.
\label{SE.6a}
\end{equation}
In particular, it follows that if the likelihood ratios $\Lambda ^1(\theta
) $ and $\Lambda ^2(\theta )$ are given on a common probability space where
they {\it coincide} as random variables, for any $\theta \in \Theta ,$ then
the experiments ${\cal E}$ and ${\cal G}$ are {\it exactly equivalent.} For
more details we refer to Nussbaum \cite{Nuss}, Proposition 2.2.

Denote by $\left( C_{[0,1]},{\cal B}(C_{[0,1]})\right) $ the measurable
space of all continuous functions on the unit interval $[0,1]$ endowed with
the uniform metric and by ${Q}_W$ the Wiener measure on $\left( C_{[0,1]},%
{\cal B}(C_{[0,1]})\right) .$

Let $P^{(i)}$, $i=1,2$ be the Gaussian shift measures on $\left( C_{[0,1]},%
{\cal B}(C_{[0,1]}),{Q}_W\right) $ induced by the following observations
\[
dX_t^{(i)}=f^{(i)}(t)dt+\frac 1{\sqrt{\sigma }}dW_t,\quad 0\leq t\leq 1,
\]
where $\sigma >0$ and $W$ is a Wiener process on $\left( C_{[0,1]},{\cal B}%
(C_{[0,1]}),{Q}_W\right) .$ Then the Hellinger distance between the measures
$P^{(1)}$ and $P^{(2)}$ satisfies the inequality (see for instance Jacod and
Shiryaev \cite{J-Sh})
\begin{equation}
H^2\left( P^{(1)},P^{(2)}\right) \leq \frac 18\sigma \int_0^1\left(
f^{(1)}(t)-f^{(2)}(t)\right) ^2dt.  \label{SE.7}
\end{equation}


\subsection{A Koml\'os-Major-Tusn\'ady approximation for independent r.v.'s%
\label{sec:KMT}}

Suppose that on the probability space $(\Omega ,{\cal F},P)$ we are given a
sequence of independent r.v.'s $X_1,...,X_n$ such that, for any $i=1,...,n,$
\[
EX_i=0
\]
and
\[
C_{\min }\leq EX_i^2\leq C_{\max },
\]
for some constants $0<C_{\min }<C_{\max }<\infty .$ Assume also that the
following Cram\'{e}r condition
\[
E\exp \{C_0|X_i|\}\leq C_1
\]
holds, for $i=1,...,n,$ with some constants $C_0>0$ and $1<C_1<\infty .$
Along with this consider that on another probability space $(\widetilde{%
\Omega },\widetilde{{\cal F}},\widetilde{P})$ we are given a sequence of
independent normal r.v.'s $N_1,...,N_n,$ with
\[
\widetilde{E}N_i=0,\quad \widetilde{E}N_i^2=EX_i^2,
\]
for $i=1,...,n.$ Let ${\cal H}(\frac 12,L)$ be the H\"{o}lder ball with
exponent $\frac 12,$ i.e.\ the set of all real valued functions $f$ defined
on the unit interval $[0,1]$ and satisfying the following conditions
\[
|f(x)-f(y)|\leq L|x-y|^{1/2},
\]
where $L>0$ and
\[
||f||_\infty \leq L/2.
\]
Let $t_i=\frac in,$ $i=1,...,n,$ be a uniform grid on the interval $[0,1].$

The following theorem is crucial for our results. The proof can be found in
 Grama and Nussbaum \cite{Gr-Nu}.

\begin{theorem}
\label{T KMT1} \label{KMT.main-theorem} Let $n\geq 2.$ A sequence
of independent r.v.'s
$\widetilde{X}_1,...,\widetilde{X}_n$ can be constructed on the probability
space $(\widetilde{\Omega },\widetilde{{\cal F}},\widetilde{P})$ such that
$\widetilde{X}_i\stackrel{d}{=}X_i,$ $i=1,...,n,$ and such that,
for $S_n(f)$ defined by
\[
S_n(f)=\sum_{i=1}^nf(t_i)(\widetilde{X}_i-N_i),
\]
we have

\[
\sup_{f\in {\cal H}(\frac 12,L)}\widetilde{P}(|S_n(f)|>x\log ^2n)\le c_1\exp
\{-c_2x\},\quad x\geq 0,
\]
where and $c_1,$ $c_2$ are constants depending only on $C_{\min },$ $C_{\max
},$ $C_0,$ $C_1,$ $L.$
\end{theorem}

\begin{remark}
It should be pointed out that in the above theorem the r.v.'s $X_i,$ $%
i=1,...,n$ are not supposed to be identically distributed nor to have
smooth distributions. To perform the construction we make use of the r.v.'s $%
N_1,...,N_n$ only, so that no additional assumptions on the probability
space $(\widetilde{\Omega },\widetilde{{\cal F}},\widetilde{P})$ are
required.
\end{remark}

Note that, under these circumstances, the construction of the r.v.'s $%
\widetilde{X}_i,$ $i=1,...,n,$ in Theorem \ref{T KMT1} appears extremely
difficult. That is why our construction of the asymptotically equivalent
Gaussian experiments can be viewed as an existence theorem rather than a
prescription for transforming the initial sample into something which
can be "assumed Gaussian" for purposes of inference.  One can expect that
a simpler construction can be performed in case of a stronger smoothness
assumption on the parameter $f(t)\in \Sigma $ and/or on the distributions of
the observed data, i.~e.~restricting somewhat the class of functions $\Sigma
$ and/or the class of allowed distributions. It remains open whether this could lead to
a constructive recipe, such as  plugging in some transform of the original data
 into an optimal estimator in the
accompanying Gaussian model and use it as an estimator for $f(t)$  in the
initial model.

A related approach is to use the Hungarian construction for constructive
purposes via closeness of sample paths of the processes, along with some
"continuity" properties of estimators. As examples for the case of kernel
density estimators we mention Rio \cite{Rio94} or Einmahl and Mason \cite
{Einm-Mas}. The precise relation to asymptotic equivalence theory remains to
be investigated.

For more details on related subjects as well as for various versions of KMT
results we refer the reader to Koml\'{o}s et al.\ \cite{KMT1}, \cite{KMT2},
Cs\"{o}rg{\H o} and R\'{e}v\'{e}sz \cite{Cs-Rev}, Sakhanenko \cite{Sakh},
Einmahl \cite{Einm}, Massart \cite{Massart}, Rio \cite{Rio93a}, \cite{Rio93b}%
, \cite{Rio94}, Koltchinskii \cite{Kol}, Einmahl and Mason \cite{Einm-Mas},
Grama and Nussbaum \cite{Gr-Nu} and to the references therein.


\section{Main results\label{sec:FormRes}}

\setcounter{equation}{0}


\subsection{Notations and formulation of the problem\label{sec:NFP}}

Let $\left( X,{\cal B}(X),\mu \right) $ a measurable space equipped with the
$\sigma $-finite measure $\mu (dx),$ where $X$ is a subset of the real line $%
R$ and ${\cal B}(X)$ is the Borel $\sigma $-field on $X.$ Let $\Theta $ be
an open (possibly infinite) interval in $R$ and ${\cal P}=\left\{ P_\theta
:\theta \in \Theta \right\} $ be an exponential family of distributions on $%
(X,{\cal B}(X),\mu )$ with parameter set $\Theta .$ The corresponding
exponential experiment we denote by ${\cal E}=\left( X,{\cal B}(X),\left\{
P_\theta :\theta \in \Theta \right\} \right) .$ Assume that conditions (\ref
{BKG.Ass 0}) and (\ref{BKG.Ass 1}) hold true. Let us introduce some further
notations.

Let ${\cal L}(\Theta _0)$ be the set of all functions $f$ defined on the
unit interval $[0,1]$ with values in the parameter set $\Theta _0$
\[
{\cal L}={\cal L}(\Theta _0)=\{f:[0,1]\rightarrow \Theta _0\},
\]
where $\Theta _0$ is the interval introduced in (\ref{BKG.Ass 1}). Let $%
{\cal H}(\beta ,L)$ be a H\"{o}lder ball, i.e. the set of functions $%
f:[0,1]\rightarrow R,$ which satisfy the conditions
\[
|f(x)|\leq L,\quad |f^{(m)}(x)-f^{(m)}(x)|\leq L|x-y|^\alpha ,
\]
for $x,y\in [0,1],$ where $\beta =m+\alpha ,$ $0<\alpha \leq 1.$ Later we
shall require that $\beta \geq \frac 12.$ Consider also the following set of
functions:
\[
\Sigma =\Sigma (\beta ,L)={\cal L}(\Theta_0 )\cap {\cal H}(\beta ,L);
\]
this will be the basic parameter set in our model.

Let $X(\theta )$ stand for a r.v.\ whose distribution is in the exponential
family ${\cal P},$ with parameter $\theta \in \Theta .$ On the unit interval
$[0,1]$ consider the time points $t_i=\frac in,$ $i=1,...,n.$ Assume that we
observe the sequence of independent r.v.'s
\begin{equation}
X_i=X(\theta _i),\quad i=1,..,n,  \label{F.4}
\end{equation}
with $\theta _i=f(t_i),$ where the "unknown" function $f$ is in the set $%
\Sigma .$ We shall prove that the statistical experiment generated by these
observations is asymptotically equivalent to an experiment of observing the
function $f$ in white noise.

Let us give another formal definition of the statistical experiments related
to the observed data $X_i,$ $i=1,...,n.$ To each time $t_i$ we associate an
exponential experiment ${\cal E}_{t_i}$ indexed by functions $f\in \Sigma $
as follows:
\[
{\cal E}_{t_i}=\left( X,{\cal B}\left( X\right) ,\left\{ P_{f(t_i)}:f\in
\Sigma \right\} \right) .
\]
Define ${\cal E}^n$ to be the product experiment ${\cal E}^n={\cal E}%
_{t_1}\otimes ...\otimes {\cal E}_{t_n}.$ In other words the experiment
corresponding to the sequence of observations $X_i,$ $i=1,...,n,$ defined by
(\ref{F.4}) is
\begin{equation}
{\cal E}^n={\cal E}(X^n)=(X^n,{\cal B}(X^n),\{P_f^n:f\in \Sigma \}),
\label{F.5}
\end{equation}
where $P_f^n$ is the product measure
\begin{equation}
P_f^n=P_{f(t_1)}\otimes ...\otimes P_{f(t_n)},\quad f\in {\cal L}.
\label{F.6}
\end{equation}
The experiment defined by (\ref{F.5}) [or equivalently by (\ref{F.4})] will
be also called {\it global} to distinguish it from the {\it local }
experiment to be introduced now.

For any fixed function $f_0$ in the parameter set $\Sigma ,$ define a
neighborhood by
\[
\Sigma _{f_0}(\gamma _n)=\{f\in \Sigma :||f-f_0||_\infty \leq \gamma _n\},
\]
where $\gamma _n\rightarrow 0$ as $n\rightarrow \infty .$ In accordance with
rate of convergence results in nonparametric statistics, the shrinking rate $%
\gamma _n$ of the neighborhood $\Sigma _{f_0}(\gamma _n)$ should be slower
than $n^{-\frac 12}.$ We will study the case where
\begin{equation}
\gamma _n=\kappa _0(n/\log n)^{-\frac \beta {2\beta +1}},  \label{F.8}
\end{equation}
where $\beta $ is the exponent in the H\"{o}lder ball ${\cal H}(\beta ,L)$
and $\kappa _0=\kappa _0(\beta )$ is some constant depending on $\beta .$
This choice can be explained in the following way. The neighborhood $\Sigma
_{f_0}(\gamma _n)$ with shrinking rate $\gamma _n$ given by (\ref{F.8})  is
such that in the experiment ${\cal E}^n$ there exists a preliminary
estimator $\widehat{f}_n$ satisfying
\begin{equation}
\sup_{f_0\in {\cal H}(\beta ,L)}P_{f_0}(\widehat{f}_n\in \Sigma
_{f_0}(\gamma _n))\rightarrow 1,\quad n\rightarrow \infty .  \label{F.8a}
\end{equation}
This property will be of use later when we globalize our local results.

The local experiment, which we will denote by ${\cal E}_{f_0}^n,$ is defined
as
\begin{equation}
{\cal E}_{f_0}^n=\left( X^n,{\cal B}(X^n),\left\{ P_f^n:f\in \Sigma
_{f_0}(\gamma _n)\right\} \right) .  \label{F.9}
\end{equation}
Let us remark that generally speaking for the sequence of nonparametric
global experiments ${\cal E}^n$ there is no Gaussian limit experiment in the
usual weak sense, since the corresponding likelihood ratios are
asymptotically degenerate. Instead it is appropriate to consider a sequence
of {\it accompanying} Gaussian experiments ${\cal G}^n,$ which is
asymptotically equivalent to the initial sequence ${\cal E}^n:$
\[
\Delta ({\cal E}^n,{\cal G}^n)\rightarrow 0,\quad n\rightarrow \infty .
\]
The same applies to the local experiments ${\cal E}_{f_0}^n.$

The corresponding accompanying Gaussian experiments will be introduced in
subsequent sections as the results are formulated. We now describe the
likelihood ratios for the experiments ${\cal E}^n$ and ${\cal E}_{f_0}^n.$
Note that, since the measure $P_\theta $ is absolutely continuous w.r.t. the
measure $\mu (dx),$
\[
\frac{P_\theta (dx)}{\mu (dx)}=\exp \{\theta U(x)-V(\theta )\},
\]
thus, for any $f\in \Sigma ,$
\[
\frac{dP_f^n}{d\mu ^n}=\prod_{i=1}^n\exp \{f(t_i)U(X_i)-V(f(t_i))\}.
\]
From this we derive that, for any $f,f_0\in \Sigma ,$ the likelihood ratio
for the measures $P_f^n$ and $P_{f_0}^n,$ corresponding to the experiment ${\cal E%
}^n={\cal E}(X^n),$ has the form
\begin{equation}
\frac{dP_f^n}{dP_{f_0}^n}=\exp
\{\sum_{i=0}^n[f(t_i)-f_0(t_i)]U(X_i)-\sum_{i=0}^n[V(f(t_i))-V(f_0(t_i))]\}.
\label{F.13}
\end{equation}
For the local experiment ${\cal E}_{f_0}^n$ the likelihood ratio has the
same form (\ref{F.13}) but with $f\in \Sigma _{f_0}(\gamma _n).$


\subsection{Local experiments: nonparametric neighborhoods\label%
{sec:LocalResults}}

We start with the local framework, since our global results are essentially
based upon the results for local experiments. For this let $f_0\in \Sigma $
be fixed. The corresponding local Gaussian experiment ${\cal G}_{f_0}^n$ is
generated by the following Gaussian observations in continuous time
\begin{equation}
dY_t^n=f(t)dt+\frac 1{\sqrt{n}}I\left( f_0(t)\right) ^{-1/2}dW_t,\quad t\in
[0,1],\quad f\in \Sigma _{f_0}(\gamma _n),  \label{L.1}
\end{equation}
where $W$ is the standard Wiener process on the probability space $\left(
C_{[0,1]},{\cal B}(C_{[0,1]}),{Q}_W\right) $. Recall that $I(\theta )$ is
the Fisher information in the local exponential family ${\cal P}$ (see
Section \ref{sec:EXP-FAM}). Denote by $Q_{f_0,f}^n$ the Gaussian shift
measure on $\left( C_{[0,1]},{\cal B}(C_{[0,1]})\right) $ induced by the
observations $(Y_t^n)_{0\leq t\leq 1}$ determined by (\ref{L.1}). Then $%
{\cal G}_{f_0}^n$ can be defined as
\begin{equation}
{\cal G}_{f_0}^n=\left( C_{[0,1]},{\cal B}(C_{[0,1]}),\left\{
Q_{f_0,f}^n:f\in \Sigma _{f_0}(\gamma _n)\right\} \right) .  \label{L.2}
\end{equation}

\begin{theorem}
\label{T L1}Assume that $\beta >\frac 12.$ Then the experiments ${\cal E}%
_{f_0}^n$ and ${\cal G}_{f_0}^n$ are asymptotically equivalent uniformly
over $f_0$ in $\Sigma $: %
\[
\sup_{f_0\in \Sigma }\Delta ({\cal E}_{f_0}^n,{\cal G}_{f_0}^n)\rightarrow
0,\quad n\rightarrow \infty .
\]
Moreover
\[
\sup_{f_0\in \Sigma }\Delta ^2({\cal E}_{f_0}^n,{\cal G}_{f_0}^n)=O\left(
n^{-\frac{2\beta -1}{2\beta +1}}(\log n)^{\frac{14\beta +5}{2\beta +1}%
}\right) .
\]
\end{theorem}

The requirement $\beta > \frac 12$ is necessary in view
of the counterexample in Brown and Low \cite{Br-Low}, remark 4.6 for the case of a Gaussian shift
exponential family {\cal P}.  We also refer the reader to other asymptotic nonequivalence results by
Efromovich and Samarov \cite{Efr-Sam} and Brown and Zhang \cite{Br-Zhang}.

We will now present a discrete version of the asymptotically equivalent
Gaussian experiment. The corresponding local experiment ${\cal Y}_{f_0}^n$
is generated by the Gaussian observations in discrete time
\begin{equation}
Y_i=f(t_i)+I\left( f_0(t_i)\right) ^{-\frac 12}\varepsilon _i,\quad
i=1,...,n,  \label{L.3}
\end{equation}
where $f\in \Sigma _{f_0}(\gamma _n)$ and $\varepsilon _i$ are i.i.d.
standard normal r.v.'s. If we denote by $G_{f_0(t_i),f(t_i)}$ the Gaussian
measure corresponding to one observation $Y_i$ of the form (\ref{L.3}),
i.e.\ the Gaussian measure on real line with mean $f(t_i)$ and variance $%
I\left( f_0(t_i)\right) ^{-1},$ then ${\cal Y}_{f_0}^n$ can be defined as
\[
{\cal Y}_{f_0}^n=\left( R^n,{\cal B}(R^n),\left\{ G_{f_0,f}^n:f\in \Sigma
_{f_0}(\gamma _n)\right\} \right) ,
\]
where
\begin{equation}
G_{f_0,f}^n=G_{f_0(t_1),f(t_1)}\otimes ...\otimes G_{f_0(t_n),f(t_n)},\quad
f\in \Sigma _{f_0}(\gamma _n).  \label{L.3a}
\end{equation}

\begin{theorem}
\label{T L1d}Assume that $\beta >\frac 12.$ Then the experiments ${\cal E}%
_{f_0}^n$ and ${\cal Y}_{f_0}^n$ are asymptotically equivalent uniformly
over $f_0$ in $\Sigma .$ Moreover
\[
\sup_{f_0\in \Sigma }\Delta ^2\left( {\cal E}_{f_0}^n,{\cal Y}%
_{f_0}^n\right) =O\left( n^{-\frac{2\beta -1}{2\beta +1}}(\log n)^{\frac{%
14\beta +5}{2\beta +1}}\right) .
\]
\end{theorem}

It is easy to see that Theorem \ref{T L1} is a consequence of Theorem \ref{T
L1d} in view of results of Brown and Low \cite{Br-Low}. Although the rate
argument is not developed there, it can easily be made explicit.


\subsection{Variance-stable form of the local approximation\label{sec:VS}}

Generally speaking there are no reasons to assume that the center of the
neighborhood $\Sigma _{f_0}(\gamma _n)$ is known to a statistician doing
nonparametric inference. That is why sometimes we would prefer to have
another form of the asymptotically equivalent Gaussian experiment ${\cal G}%
_{f_0}^n,$ in which the expression $I(f_0(t))$ does not appear. It turns out
that such a form of the Gaussian accompanying experiment does exist. We will
call it {\it variance-stable form}, since it involves the
variance-stabilizing transformation pertaining to the exponential family.
This form of the asymptotically equivalent local Gaussian experiment can be
introduced by means of a transformation of the parameter space given by any
function $\Gamma(\theta):\Theta \rightarrow R $ satisfying
\begin{equation}
\frac{d }{d\theta }\Gamma(\theta)=\sqrt{I(\theta)},  \label{VS.001}
\end{equation}
where $I(\theta)$ is the Fisher information in the local exponential family $%
{\cal P} $ [see (\ref{BKG.4}) in Section \ref{sec:EXP-FAM}]. To show its
connection with the variance-stabilizing transformation introduced in
Section \ref{sec:VST} we need some notations.
Let, as before, $V(\theta )$ be
the cumulant generating function associated with the exponential family $%
{\cal P}. $ Set, for brevity, $b(\theta )=V^{\prime }(\theta ),$ $\theta \in
\Theta .$ It follows from the assumption (\ref{BKG.Ass 0}) that $b(\theta )$
is an increasing differentiable function on the open interval $\Theta .$
Denote by $\Lambda $ the range of $b(\theta ),$ i.e.\ $\Lambda =\left\{
b(\theta ):\theta \in \Theta \right\} .$ It is clear that $\Lambda $ is also
an open interval in $R.$ Let $a(\lambda ),$ $\lambda \in \Lambda $ be the
inverse of $b(\theta ),$ $\theta \in \Theta ,$ i.e.\
\begin{eqnarray}
a(\lambda )=\inf \left\{ \theta \in \Theta :b(\theta )>\lambda \right\}
,\quad \lambda \in \Lambda ,  \nonumber
\end{eqnarray}
which obviously is an increasing differentiable function on $\Lambda .$
Another equivalent way to define $a(\lambda )$ would be to put $a(\lambda
)=T^{\prime }(\lambda ),$ where $T(\lambda )$ is the Legendre transformation
of the function $V(\theta ):$
\begin{eqnarray}
T(\lambda )=\inf \left\{ \lambda \theta -V(\theta ):\lambda \in \Lambda
\right\} .  \nonumber
\end{eqnarray}
It is also easy to see that $a(\lambda )$ satisfies the equation
\begin{equation}
a^{\prime }(\lambda )=I(a(\lambda ))^{-1},\quad \lambda \in \Lambda .
\label{VS.3a}
\end{equation}
Let $F(\lambda )$ be any function on $\Lambda ,$ having the property
\begin{equation}
F^{\prime }(\lambda )=\sqrt{a^{\prime }(\lambda )},\quad \lambda \in \Lambda.
\label{VS.4}
\end{equation}
The relations (\ref{VS.3a}) and (\ref{VS.4}) show that $F(\lambda ),$ $%
\lambda \in \Lambda $ coincides with the variance-stabilizing transformation
defined in Section \ref{sec:VST}. The functions $b(\theta ):\Theta
\rightarrow \Lambda $ and $F(\lambda ):\Lambda \rightarrow R$ define the
transformation $\Gamma (\theta )=F(b(\theta )), $ which can be easily seen
to satisfy (\ref{VS.001}), in view of (\ref{VS.3a}) and (\ref{VS.4}).

Let $\widehat{{\cal G}}_{f_0}^n$ be the Gaussian experiment
generated by the observations
\begin{equation}
d\widehat{Y}_t^n=\Gamma \left( f(t)\right) dt+\frac 1{\sqrt{n}}dW_t,\quad
t\in [0,1],\quad f\in \Sigma _{f_0}(\gamma _n),  \label{VS.5}
\end{equation}
i.e.\
\begin{eqnarray}
\widehat{{\cal G}}_{f_0}^n=\left( C_{[0,1]},{\cal B}(C_{[0,1]}),\left\{
\widehat{Q}_f^n:f\in \Sigma _{f_0}(\gamma _n)\right\} \right) ,  \nonumber
\end{eqnarray}
where $\widehat{Q}_f^n$ is the Gaussian shift measure on $\left( C_{[0,1]},%
{\cal B}(C_{[0,1]})\right) $ induced by the observations $(\widehat{Y}%
_t^n)_{0\leq t\leq 1}$ determined by (\ref{VS.5}). Let ${\cal G}_{f_0}^n$ be
the Gaussian experiment defined in (\ref{L.1}) and (\ref{L.2}).

\begin{theorem}
\label{T VS1}Let $\beta >\frac 12.$ Then the experiments ${\cal G}_{f_0}^n$
and $\widehat{{\cal G}}_{f_0}^n$ are asymptotically equivalent uniformly in $%
f_0\in \Sigma .$ Moreover the Le Cam distance between ${\cal G}_{f_0}^n$ and
$\widehat{{\cal G}}_{f_0}^n$ satisfies
\[
\sup_{f_0\in \Sigma }\Delta ^2\left( {\cal G}_{f_0}^n,\widehat{{\cal G}}%
_{f_0}^n\right) =O\left( n^{-\frac{2\beta -1}{2\beta +1}}\left( \log
n\right) ^{\frac{4\beta }{2\beta +1}}\right) .
\]
\end{theorem}

As an immediate consequence of Theorems \ref{T L1} and \ref{T VS1} we get
the following result.

\begin{theorem}
\label{T VS2}Let $\beta >\frac 12.$ Then the experiments ${\cal E}_{f_0}^n$
and $\widehat{{\cal G}}_{f_0}^n$ are asymptotically equivalent uniformly in $%
f_0\in \Sigma $ and the Le Cam distance between ${\cal E}_{f_0}^n$ and $%
\widehat{{\cal G}}_{f_0}^n$ satisfies
\[
\sup_{f_0\in \Sigma }\Delta ^2\left( {\cal E}_{f_0}^n,\widehat{{\cal G}}%
_{f_0}^n\right) =O\left( n^{-\frac{2\beta -1}{2\beta +1}}\left( \log
n\right) ^{\frac{14\beta +5}{2\beta +1}}\right) .
\]
\end{theorem}

We turn to the discrete version of the asymptotically equivalent Gaussian
experiment. The local experiment $\widehat{{\cal Y}}_{f_0}^{n}$ is generated
by the Gaussian observations in discrete time
\begin{equation}
\widehat{Y}_i=\Gamma \left( f(t_i)\right) +\varepsilon _i,\quad f\in \Sigma
_{f_0}(\gamma _n),\quad i=1,...,n,  \label{VS.6}
\end{equation}
where $\varepsilon _i$ are i.i.d.\ standard normal r.v.'s. If we denote by $%
G_{f(t_i)}$ the Gaussian measure corresponding to one observation $Y_i$ of
the form (\ref{VS.6}), i.e.\ the Gaussian measure on the real line with mean
$\Gamma \left( f(t_i)\right) $ and variance $1,$ then $\widehat{{\cal Y}}%
_{f_0}^{n}$ can be defined as
\[
\widehat{{\cal Y}}_{f_0}^{n}= \left( R^n,{\cal B}(R^n), \left\{ \widehat
G_f^{n}:f\in \Sigma _{f_0}(\gamma _n)\right\} \right) ,
\]
where
\[
\widehat G_f^{n}= G_{f(t_1)}\otimes ...\otimes G_{f(t_n)},\quad f\in \Sigma
_{f_0}(\gamma _n).
\]

\begin{theorem}
\label{T VS2d}Assume that $\beta >\frac 12.$ Then the experiments ${\cal E}%
_{f_0}^n$ and $\widehat{{\cal Y}}_{f_0}^n$ are asymptotically equivalent
uniformly over $f_0$ in $\Sigma .$ Moreover
\[
\sup_{f_0\in \Sigma }\Delta ^2\left( {\cal E}_{f_0}^n,\widehat{{\cal Y}}%
_{f_0}^n\right) =O\left( n^{-\frac{2\beta -1}{2\beta +1}}(\log n)^{\frac{%
14\beta +5}{2\beta +1}}\right) .
\]
\end{theorem}

In the above theorems the initial exponential experiment ${\cal E},$ which
generates ${\cal E}_{f_0}^n,$ is assumed to be in its canonical form. The
variance-stable Gaussian approximation appears in an equivalent but a little
more pleasant form (as we will see in Section \ref{sec:Examples}), if the
experiment ${\cal E}$ is {\it naturally parametrized,} i.e.\ if ${\cal E},$
given in its canonical form by (\ref{BKG.1}), is reparametrized by means of
the one-to-one map $b(\theta ):\Theta _0\rightarrow \Lambda _0.$
Introduce the set of functions
\[
\overline{\Sigma }=\left\{ g=b\circ f:f\in \Sigma \right\}
\]
and, for any $g_0\in \overline{\Sigma },$ the neighborhoods
\[
\overline{\Sigma }_{g_0}(\gamma _n)=\left\{ g\in \overline{\Sigma }:\left\|
g-g_0\right\| _\infty \leq \overline{c}_0\gamma _n\right\} ,
\]
with some constant $\overline{c}_0$ depending $I_{\max },$ $\varepsilon _0.$
Let $\overline{{\cal E}}_{g_0}^n$ be the corresponding nonparametric product
experiment, defined analogously to (\ref{F.9}) and (\ref{F.6}):
\[
\overline{{\cal E}}_{g_0}^n=\left( X^n,{\cal B}(X^n),\left\{ \overline{P}%
_g^n:g\in \overline{\Sigma }_{g_0}(\gamma _n)\right\} \right) ,
\]
with $\overline{P}_g^n=P_{f}^n$ for  $f=a\circ g$, where $P_f^n$ is
the product measure defined by (\ref{F.6}). The accompanying Gaussian
experiment $\overline{{\cal G}}_{g_0}^n$ is defined by the observations
\[
d\overline{Y}_t^n=F(g(t))dt + \frac 1{\sqrt{n}}dW_t,\quad t\in [0,1],\quad g\in
\overline{\Sigma }_{g_0}(\gamma _n).
\]

\begin{theorem}
\label{T VS3}Let $\beta >\frac 12.$ Then the experiments $\overline{{\cal E}}%
_{g_0}^n$ and $\overline{{\cal G}}_{g_0}^n$ are asymptotically equivalent
uniformly in $g_0\in \overline{\Sigma }$ and the Le Cam distance between $%
\overline{{\cal E}}_{g_0}^n$ and $\overline{{\cal G}}_{g_0}^n$ satisfies
\[
\sup_{g_0\in \overline{\Sigma }}\Delta ^2\left( \overline{{\cal E}}_{g_0}^n,%
\overline{{\cal G}}_{g_0}^n\right) =O\left( n^{-\frac{2\beta -1}{2\beta +1}%
}\left( \log n\right) ^{\frac{14\beta +5}{2\beta +1}}\right) .
\]
\end{theorem}

The discrete version of this theorem is straightforward.

\subsection{Global experiments\label{sec:GlobalResults}}

The variance-stable form of the accompanying {\it local} Gaussian
experiments allows to construct an accompanying {\it global} Gaussian
experiment. The main idea is to substitute a preliminary estimator
satisfying (\ref{F.8a}) for the unknown function $f_0,$ around which the
local experiment is built. Such an estimator is provided by Lemma \ref{L ES1}
(see Section \ref{sec:PRES}). In the variance-stable form of the local
experiment given by (\ref{VS.5}) (or (\ref{VS.6}) in the discrete case) the
unknown function $f_0$ does not show up in the distributions themselves but
appears only as a center of the parametric neighborhood. This will imply in
the sequel that the globalized experiment does not depend on the specific
preliminary estimator used; thus a convenient closed form global
approximation for the original regression experiment ${\cal E}^n$ can be
obtained.

Let $\Gamma (\theta ):\Theta _0\rightarrow R$ be the transformation given by
(\ref{VS.001}). Let the global Gaussian experiment ${\cal G}^n$ be defined
as
\[
{\cal G}^n=\left( R^n,{\cal B}(R^n),\left\{ Q_f^n:f\in \Sigma \right\}
\right) ,
\]
where $Q_f^n$ is the Gaussian shift measure induced by the observations
\[
dY_t^n=\Gamma \left( f(t)\right) dt+\frac 1{\sqrt{n}}dW_t,\quad f\in \Sigma
,\quad t\in [0,1].
\]

\begin{theorem}
\label{T G2}Let $\beta >\frac 12.$ Then the global experiments ${\cal E}^n$
and ${\cal G}^n$ are asymptotically equivalent: $\Delta ({\cal E}^n,{\cal G}%
^n)\rightarrow 0$ as $n\rightarrow \infty .$ Moreover for $\beta \in (\frac
12,1)$
\[
\Delta ^2({\cal E}^n,{\cal G}^n)=O\left( n^{-\frac{2\beta -1}{2\beta +1}%
}(\log n)^{\frac{14\beta +5}{2\beta +1}}\right) .
\]
\end{theorem}

As a particular case this theorem gives us the main result in Brown and Low
\cite{Br-Low}. We will discuss the case of normal observations and other
examples of interest in the next section.

We present also a discrete version of the asymptotically equivalent global
Gaussian experiment in its variance-stable form. The Gaussian experiment $%
{\cal Y}^n$ is generated by the observations in discrete time
\begin{equation}
Y_i=\Gamma \left( f(t_i)\right) +\varepsilon _i,\quad f\in \Sigma ,\quad
i=1,...,n,  \label{G.8}
\end{equation}
where $\varepsilon _i$ are i.i.d.\ standard normal r.v.'s. If we denote by $%
G_{f(t_i)}$ the Gaussian measure corresponding to one observation $Y_i$ of
the form (\ref{G.8}), $i=1,...,n,$ i.e. the Gaussian measure on the real
line with mean $f(t_i)$ and variance $1,$ then ${\cal Y}^n$ can be defined
as
\begin{equation}
{\cal Y}^n=\left( R^n,{\cal B}(R^n),\left\{ G_f^n:f\in \Sigma \right\}
\right) ,  \label{G.9}
\end{equation}
where
\begin{equation}
G_f^n=G_{f(t_1)}\otimes ...\otimes G_{f(t_n)},\quad f\in \Sigma .
\label{G.10}
\end{equation}

\begin{theorem}
\label{T G2d}Let $\beta >\frac 12.$ Then the global experiments ${\cal E}^n$
and ${\cal Y}^n$ are asymptotically equivalent: $\Delta ({\cal E}^n,{\cal Y}%
^n)\rightarrow 0$ as $n\rightarrow \infty .$ Moreover for $\beta \in (\frac
12,1)$
\[
\Delta ^2({\cal E}^n,{\cal Y}^n)=O\left( n^{-\frac{2\beta -1}{2\beta +1}%
}(\log n)^{\frac{14\beta +5}{2\beta +1}}\right) .
\]
\end{theorem}

The case when ${\cal E}^{n}$ is naturally parametrized is similar to Theorem
\ref{T VS3}. The proofs for the passage from local to global are in Section \ref{sec:GLOB}%
.

\subsection{Local experiments: almost $n^{-1/2}$-neighborhoods \label%
{sec:almostpar}}

It turns out that the key point in the study of asymptotic equivalence of
local experiments is the behavior of its fragments over shrinking time
intervals of length
\[
\delta _n=\gamma _n^{1/\beta }=\kappa _0^{1/\beta }(n/\log n)^{-\frac
1{2\beta +1}}.
\]
After a rescaling we arrive at local experiments parametrized by functions $%
f $ of the form $f=f_0+\gamma _n^{*}g,$ where $g\in \Sigma $ and
\begin{equation}
\gamma _n^{*}=\kappa _0^{*}(n/\log n)^{-\frac 12},  \label{PN.2}
\end{equation}
with some $\kappa _0^{*}>0.$ We will present the corresponding results,
since they are of independent interest. Before stating these results we
introduce the necessary notations.

Let $f_0\in {\cal L}$ (see Section \ref{sec:NFP}). In the sequel all
functions $g\in \Sigma $ are assumed to be such that $f=f_0+\gamma _n^{*}g\in
{\cal L}.$ For any $g\in \Sigma $ set
\[
P_{f_0,g}^n=P_{f_0+\gamma _n^{*}g}^n,
\]
where $P_f^n$ is the product measure defined in (\ref{F.6}) and consider the
local experiment
\begin{equation}
{\cal E}_{f_0}^{*,n}=\left( X^n,{\cal B}(X^n),\left\{ P_{f_0,g}^n:g\in
\Sigma \right\} \right) .  \label{PN.2a}
\end{equation}
In the same way we introduce the accompanying sequence of Gaussian
experiments:
\[
{\cal G}_{f_0}^{*,n}=\left( C_{[0,1]}^n,{\cal B}(C_{[0,1]}^n),\left\{
Q_{f_0,g}^{*,n}:g\in \Sigma \right\} \right) ,
\]
where $Q_{f_0,g}^{*,n}=Q_{f_0,f_0+\gamma _n^{*}g}^n$ and $Q_{f_0,f}^n$ is
generated by the Gaussian observations (\ref{L.1}), with $f_0\in {\cal L}$
and $f=f_0+\gamma _n^{*}g,$ $g\in \Sigma .$

\begin{theorem}
\label{T LD1}Assume that $\beta >\frac 12.$ Then the experiments ${\cal E}%
_{f_0}^{*,n}$ and ${\cal G}_{f_0}^{*,n}$ are asymptotically equivalent
uniformly over $f_0$ in ${\cal L}.$ Moreover
\[
\sup_{f_0\in {\cal L}}\Delta ^2\left( {\cal E}_{f_0}^{*,n},{\cal G}%
_{f_0}^{*,n}\right) =O\left( n^{-1}(\log n)^7\right) .
\]
\end{theorem}

Let ${\cal Y}_{f_0}^{*,n}$ be the local experiment corresponding to the
independent Gaussian observations in discrete time $Y_i,$ $i=1,...,n,$
defined by (\ref{L.3}). More precisely, let
\[
G_{f_0,g}^{*,n}=G_{f_0,f_0+\gamma _n^{*}g}^n,
\]
where the probability measure $G_{f_0,f}^n$ is generated by the Gaussian
observations (\ref{L.3}), with $f_0\in {\cal L}$ and $f=f_0+\gamma _n^{*}g,$
$g\in \Sigma ,$ and set
\begin{equation}
{\cal Y}_{f_0}^{*,n}=\left( R^n,{\cal B}(R^n),\left\{ G_{f_0,g}^{*,n}:g\in
\Sigma \right\} \right) .  \label{PN.2b}
\end{equation}

\begin{theorem}
\label{T LD1d}Assume that $\beta >\frac 12.$ Then the experiments ${\cal E}%
_{f_0}^{*,n}$ and ${\cal Y}_{f_0}^{*,n}$ are asymptotically equivalent
uniformly over $f_0$ in ${\cal L}.$ Moreover
\[
\sup_{f_0\in {\cal L}}\Delta ^2({\cal E}_{f_0}^{*,n},{\cal Y}%
_{f_0}^{*,n})=O\left( n^{-1}(\log n)^7\right) .
\]
\end{theorem}


\section{Examples and applications\label{sec:Examples}}

\setcounter{equation}{0}


The most striking form for the asymptotically equivalent Gaussian
approximation in the following examples will be obtained  if the initial
exponential experiment ${\cal E}$ is taken under its {\it natural
parametrization.} By a natural parametrization we mean the following. Let $%
{\cal P}=\left\{ P_\theta :\theta \in \Theta \right\} $ be an exponential
family on $(X,{\cal B}(X),\mu )$ in the canonical form, i.~e.~whose
Radon-Nikodym derivatives $dP_\theta /d\mu $ are defined by (\ref{BKG.1}).
Let $V(\theta )$ be its cumulant generating function. It is clear (see also
Section \ref{sec:VS}) that $b(\theta )=V^{\prime }(\theta )$ is a one-to-one
map from $\Theta _0$ to $\Lambda _0.$ If we reparametrize the family ${\cal P%
}$ by means of the map $b(\theta ),$ we will call the family ${\cal P}$ {\it %
naturally parametrized.} Indeed the parameter $\lambda =b(\theta )$ is the
''natural'' parameter in many specific families: mean or variance for normal
distributions, intensity for exponential, gamma, or Poisson distributions,
probability of success for Bernoulli and binomial distributions
etc.\thinspace .


\subsection{Gaussian observations: unknown mean\label{sec:GO-mean}}

Let ${\cal P}=\left\{ P_\theta :\theta \in \Theta \right\} $ be the family
of normal distributions on the real line $X=R$ with mean $\theta \in \Theta
=R$ and variance $1.$ The normal distribution $P_\theta (dx)$ can be written
\[
P_\theta (dx)=e^{\theta x-V(\theta )}\mu (dx),
\]
$\mu (dx)$ being the standard normal distribution
\[
\mu (dx)=\frac 1{\sqrt{2\pi }}e^{-\frac{x^2}2}dx
\]
and $U(x)=x,$ $V(\theta )=\theta ^2/2.$ Then the corresponding Fisher
information $I(\theta )=V^{\prime \prime }(\theta )\equiv 1,$ so the
condition (\ref{BKG.Ass 1}) holds true with $\Theta _0=R,$ $I_{\max
}=I_{\min }=1,$ $\varepsilon _0>0.$ Hence the parameter set $\Sigma $
coincides with the H\"{o}lder ball ${\cal H}(\beta ,L).$

Assume that our observations $X_i=X(\theta _i),$ $i=1,...,n,$ are
independent normal with mean $\theta _i=f(t_i),$ $t_i=i/n$ and variance $1,$
where the function $f$ is in the H\"{o}lder ball $\Sigma ={\cal H}(\beta
,L). $ Let us remark that these observations correspond to the regression
model
\[
X_i=f(t_i)+\varepsilon _i,\quad i=1,...,n,
\]
with i.i.d.\ standard normal r.v.'s $\varepsilon _i,$ $i=1,...,n.$

Let ${\cal E}_{f_0}^n$ be the local experiment generated by the sequence of
observations $X_i,$ $i=1,...,n,$ with $f\in \Sigma _{f_0}(\gamma _n),$ for
some $f_0\in \Sigma .$ According to Theorem \ref{T L1}, for any $\beta
>\frac 12,$ the experiment ${\cal E}_{f_0}^n$ is asymptotically equivalent
to the local experiment ${\cal G}_{f_0}^n$ generated by the observations
\begin{equation}
dY_t^n=f(t)dt+\frac 1{\sqrt{n}}dW_t,\quad t\in [0,1],  \label{GM.2}
\end{equation}
where $f\in \Sigma _{f_0}(\gamma _n).$

The global form of the asymptotically equivalent Gaussian experiment is
given also by (\ref{GM.2}) but with $f\in \Sigma .$ Thus, we recover the
main result of Brown and Low \cite{Br-Low}.


\subsection{Gaussian observations: unknown variance\label{sec:GO-var}}

Let ${\cal P}=\left\{ P_\lambda :\lambda \in \Lambda \right\} $ be the
family of normal distributions on the real line $X=R$ with mean $0$ and
variance $\lambda \in \Lambda \equiv [0,\infty ).$ The normal distribution $%
P_\lambda (dx)$ has the form
\begin{equation}
P_\lambda (dx)=\frac 1{\sqrt{2\pi \lambda }}\exp \{-\frac{x^2}{2\lambda }%
\}dx.  \label{EGV.0}
\end{equation}
After the reparametrization $\theta =-1/\lambda ,$ $\theta \in \Theta
=(-\infty ,0],$ we obtain the linearly indexed exponential model
\begin{equation}
P_\theta (dx)=\exp \{\theta U(x)-V(\theta )\}\mu (dx),  \label{EGV.1}
\end{equation}
where $\mu (dx)$ is the Lebesgue measure on the real line and $U(x)=x^2/2,$ $%
V(\theta )=-\frac 12\log (-\frac \theta {2\pi }).$ The corresponding Fisher
information is $I(\theta )=V^{\prime \prime }(\theta )=\frac 12\theta ^{-2},$
so condition (\ref{BKG.Ass 1}) holds true with $\Theta _0=[\theta _{\min
},\theta _{\max }],$ for some constants $-\infty <\theta _{\min }<\theta
_{\max }<0$ and $\varepsilon _0$ small enough. Then the parameter set $%
\Sigma $ contains all the functions $f(t)$ from the H\"{o}lder ball ${\cal H}%
(\beta ,L),$ which satisfy $\theta _{\min }\leq f(t)\leq \theta _{\max }.$

Consider a sequence of independent normal observations
\begin{equation}
X_i=X(\theta _i),\quad i=1,...,n,  \label{EGV.2}
\end{equation}
of means zero and ''canonical variances'' $\theta _i=f(t_i),$ $t_i=i/n,$
where the unknown function $f$ is in the set $\Sigma .$ Let $f_0(t)\in
\Sigma $ and let ${\cal E}_{f_0}^n$ be the local experiment
generated by the observations (\ref{EGV.2}).
Then, by Theorem \ref{T L1}, the experiment
${\cal E}_{f_0}^n$ is asymptotically equivalent to the Gaussian experiment
${\cal G}_{f_0}^n$ generated by the observations
\[
dY_t^n=f(t)dt-\frac{\sqrt{2}}{\sqrt{n}}f_0(t)dW_t,\quad t\in [0,1],
\]
for $f\in \Sigma _{f_0}(\gamma _n),$ $W$ being the standard Wiener process.

For the variance-stable form we easily compute $b(\theta )=V^{\prime
}(\theta )=-\frac 1{2\theta }$ and $F(\lambda )=\log \lambda /\sqrt{2}.$
Thus, by Theorem \ref{T VS1}, the variance-stable accompanying Gaussian
experiment is given by
\begin{equation}
dY_t^n=\frac 1{\sqrt{2}}\log \left( -\frac 1{2f(t)}\right) dt-\frac 1{\sqrt{n%
}}dW_t,\quad t\in [0,1].  \label{EGV.4*}
\end{equation}
Note that in the above formula $f(t)$ and $f_0(t)$ are less than $\theta
_{\max }<0.$

A more compact form for the accompanying Gaussian experiments is obtained
using the natural parametrization. If we reparametrize the exponential
family (\ref{EGV.1}) by means of the map $b(\theta )=-\frac 1{2\theta },$ we
recover its original form given by (\ref{EGV.0}). Let $\overline{\Sigma }=%
\overline{\Sigma }(\beta ,L)$ be the set of all functions $g$ from the
H\"{o}lder ball ${\cal H}(\beta ,L),$ which satisfy $\lambda _{\min }\leq
g(t)\leq \lambda _{\max },$ where $0<\lambda _{\min }<\lambda _{\max
}<\infty .$ For any function $g_0\in \overline{\Sigma },$ let $\overline{%
\Sigma }_{g_0}(\gamma _n)$ be its neighborhood of radius $\gamma _n$ [see (%
\ref{F.8})] in $\overline{\Sigma }.$ Denote by $\overline{{\cal {E}}}%
_{g_0}^n $ the local experiment generated by the observations (\ref{EGV.2})
with $\theta _i=-1/g(t_i),$ $g\in \overline{\Sigma }_{g_0}(\gamma _n).$ This
experiment can be regarded as generated by the observations
\[
X_i=\sqrt{g(t_i)}\varepsilon _i,\quad i=1,...,n,
\]
with $g\in \overline{\Sigma }_{g_0}(\gamma _n),$ where $\varepsilon _i$ are
i.i.d.\ standard normal r.v.'s. Then, for any $\beta >\frac 12,$ the
experiment $\overline{{\cal E}}_{g_0}^n$ is asymptotically equivalent to the
local Gaussian experiment generated by the observations
\[
dY_t^n=-\frac 1{g(t)}dt+\frac{\sqrt{2}}{\sqrt{n}g_0(t)}dW_t,\quad 0\leq
t\leq 1,
\]
with $g\in \overline{\Sigma }_{g_0}(\gamma _n),$ $W$ being
the standard Wiener
process on $(C_{[0,1]},{\cal B}(C_{[0,1]}),{Q}_W).$

The corresponding variance-stable form is determined by the equation (cf.
Theorem \ref{T VS3})
\begin{equation}
dY_t^n=\frac 1{\sqrt{2}}\log g(t)dt+\frac 1{\sqrt{n}}dW_t,\quad t\in [0,1].
\label{EGV.6}
\end{equation}
Global variants of the accompanying Gaussian experiments are also given by (%
\ref{EGV.4*}) and (\ref{EGV.6}), with extended parameter space.


\subsection{Poisson observations\label{sec:PO}}

Let $X=\{0,1,...\}$ and $\mu (dx)$ be the $\sigma $-finite measure on $X$
assigning $1/x!$ to each point $x\in X.$ Let us consider the case when $%
{\cal P}$ is the family of Poisson distributions $P_\lambda (x),$ $x\in X,$
with intensity $\lambda \in (0,\infty ).$ After the reparametrization $%
\theta =\log \lambda ,$ $\theta \in \Theta =(-\infty ,\infty ),$ we get the
canonical form
\begin{eqnarray}
P_\theta (dx)=e^{\theta x-V(\theta )}\mu (dx),  \nonumber
\end{eqnarray}
with $V(\theta )=e^\theta .$ The corresponding Fisher information is $%
I(\theta )=e^\theta .$ It is clear that condition (\ref{BKG.Ass 1}) holds
true with $\Theta _0=[\theta _{\min },\theta _{\max }],$ where $-\infty
<\theta _{\min }<\theta _{\max }<\infty.$ Then the parameter set $\Sigma $
contains all the functions $f$ from the H\"older ball ${\cal H}(\beta, L),$
which satisfy $\theta _{\min }\leq f(t)\leq \theta _{\max }.$

Consider a sequence of independent Poisson observations
\begin{equation}
X_i=X(\theta _i),\quad i=1,...,n,  \label{EP.2}
\end{equation}
with "canonical intensity" $\theta _i=f(t_i),$ $t_i=i/n,$ where the unknown
function $f$ is assumed to be in the set $\Sigma .$ Let $f_0(t)\in \Sigma $
and ${\cal E}_{f_0}^n$ be the local experiment generated by
the observations (\ref{EP.2}) with $f\in \Sigma _{f_0}(\gamma _n).$
Then, according to Theorem
\ref{T L1}, the experiment ${\cal E}_{f_0}^n$ is asymptotically equivalent
to the local Gaussian experiment ${\cal G}_{f_0}^n$
generated by the observations
\begin{equation}
dY_t^n=f(t)dt+\frac 1{\sqrt{n}}e^{-f_0(t)/2}dW_t,\quad t\in [0,1],
\label{EP.3}
\end{equation}
where $f\in \Sigma _{f_0}(\gamma _n).$

A variance-stable form of the observations (\ref{EP.3}) can be obtained, if
we note that $b(\theta )=e^\theta $ and $F(\lambda )=2\sqrt{\lambda }.$ In
view of Theorem \ref{T VS1} the experiment ${\cal G}_{f_0}^n$ is
asymptotically equivalent to the experiment $\widehat{{\cal G}}_{f_0}^n$
given by the observations
\begin{equation}
d\widehat{Y}_t^n=2\sqrt{\log f(t)}dt+\frac 1{\sqrt{n}}dW_t,\quad t\in [0,1],
\label{EP.4}
\end{equation}
where $f\in \Sigma _{f_0}(\gamma _n).$

In terms of the original parameter $\lambda =e^\theta $ these results can be
formulated as follows. Let $\overline{\Sigma }=\overline{\Sigma }(\beta ,L)$
be the set of all functions $g$ from the H\"{o}lder ball ${\cal H}(\beta
,L), $ which satisfy $\lambda _{\min }\leq g(t)\leq \lambda _{\max },$ where
$0<\lambda _{\min }<\lambda _{\max }<\infty .$ For any function $g_0\in
\overline{\Sigma },$ let $\overline{\Sigma }_{g_0}(\gamma _n)$ be its
neighborhood of radius $\gamma _n$ [see (\ref{F.8})] in $\overline{\Sigma }.$
Denote by $\overline{{\cal E}}_{g_0}^n$ the local experiment generated by the
independent Poisson observations
\[
X_i=X(\lambda _i),\quad i=1,...,n,
\]
with unknown intensities $\lambda _i=g(t_i),$ where $g\in \overline{\Sigma }%
_{g_0}^n(\gamma _n).$ Then, for any $\beta >\frac 12,$ the experiment $%
\overline{{\cal E}}_{g_0}^n$ is asymptotically equivalent to the local
Gaussian experiment generated by the observations
\[
dY_t^n=\log g(t)dt+\frac 1{\sqrt{ng_0(t)}}dW_t,\quad t\in [0,1],
\]
where $g\in \overline{\Sigma }_{g_0}(\gamma _n),$ $W$ being the standard
Wiener process on $(C_{[0,1]},{\cal B}(C_{[0,1]}),{Q}_W).$

The variance-stable result is furnished by Theorem \ref{T VS3}: an
accompanying Gaussian experiment (local or global) is also given by the
equation
\[
dY_t^n=2\sqrt{g(t)}dt+\frac 1{\sqrt{n}}dW_t,\quad t\in [0,1].
\]


\subsection{Bernoulli observations\label{sec:BO}}

Let ${\cal P}$ be the family of Bernoulli distributions $P_\lambda (x),$ $%
x\in X=\{0,1\}$ with $P_\lambda (1)=\lambda ,$ $P_\lambda (0)=1-\lambda , $ $%
\lambda \in (0,1).$ After the reparametrization $\theta =\log \frac \lambda
{1-\lambda },$ $\theta \in \Theta =(-\infty ,\infty )$ we arrive at the
following canonical form
\[
P_\theta (x)=e^{\theta x-V(\theta )},\quad x\in X,
\]
where $V(\theta )=\log (1+e^\theta ).$ The corresponding Fisher information
is $I(\theta )=e^\theta /(1+e^\theta )^2.$ One can easily check that the
condition (\ref{BKG.Ass 1}) holds true with $\Theta _0=[\theta _{\min
},\theta _{\max }],$ where $-\infty <\theta _{\min }<\theta _{\max }<\infty
. $ Then the parameter set $\Sigma $ contains all the functions $f(t)$ from
the H\"older ball ${\cal H}(\beta ,L),$ which satisfy $\theta _{\min }\leq
f(t)\leq \theta _{\max }.$

Consider a sequence of independent Bernoulli observations
\begin{equation}
X_i=X(\theta _i),\quad i=1,...,n  \label{EB.1}
\end{equation}
with canonical parameter $\theta _i=f(t_i),$ $t_i=i/n,$ where the unknown
function $f$ is in the parameter set $\Sigma .$ Let $f_0(t)\in \Sigma $ and $%
{\cal E}_{f_0}^n$ be the local experiment generated by the observations
(\ref{EB.1}), with $f\in \Sigma _{f_0}(\gamma _n).$
Then, according to Theorem
\ref{T L1}, the experiment ${\cal E}_{f_0}^n$ is asymptotically equivalent
to the local Gaussian experiment ${\cal G}_{f_0}^n$ generated by the
observations
\[
dY_t^n=f(t)dt+\frac 1{\sqrt{n}}\frac{1+e^{f_0(t)}}{e^{f_0(t)/2}}dW_t,\quad
t\in [0,1],
\]
with $f\in \Sigma _{f_0}(\gamma _n).$

For the variance-stable form we compute $b(\theta )=e^\theta /(1+e^\theta )$
and $F(\lambda )=2\arcsin \sqrt{\lambda }.$ Then, by Theorem \ref{T VS1},
the variance-stable accompanying Gaussian experiment is associated with the
equation
\[
dY_t^n=2\arcsin \sqrt{\frac{e^{f(t)}}{1+e^{f(t)}}}\,dt+\frac 1{\sqrt{n}%
}dW_t,\quad t\in [0,1].
\]

In term of the original parameter $\lambda =e^\theta /(1+e^\theta )$ this
result can be formulated in the following way. Let $\overline{\Sigma }=%
\overline{\Sigma }(\beta ,L)$ be the set of all functions $g$ from the
H\"{o}lder ball ${\cal H}(\beta ,L),$ which satisfy $\lambda _{\min }\leq
g(t)\leq \lambda _{\max },$ where $0<\lambda _{\min }<\lambda _{\max }<1.$
For any function $g_0\in \overline{\Sigma },$ let $\overline{\Sigma }%
_{g_0}(\gamma _n)$ be its neighborhood of radius $\gamma _n$ [see (\ref{F.8}%
)] in $\overline{\Sigma }.$ Denote by $\overline{{\cal E}}_{g_0}^n$ the
local experiment generated by the observations (\ref{EB.1}),
with $\theta _i=\log\frac{g(t_i)}{1-g(t_i)},$
$g\in \overline{\Sigma }_{g_0}(\gamma _n).$ Then,
for any $\beta >\frac 12,$ the experiment $\overline{{\cal E}}_{g_0}^n$ is
asymptotically equivalent to the local Gaussian experiment generated by the
observations
\[
dY_t^n=\log \frac{g(t)}{1-g(t)}dt+\frac 1{\sqrt{ng_0(t)(1-g_0(t))}%
}dW_t,\quad t\in [0,1],
\]
where $g\in \overline{\Sigma }_{g_0}(\gamma _n),$ $W$ being the standard
Wiener process on $(C_{[0,1]},{\cal B}(C_{[0,1]}),{Q}_W).$

The variance-stable form of the Gaussian accompanying experiment
parametrized by $g$ is (according to Theorem \ref{T VS3})
\[
dY_t^n=2\arcsin \sqrt{g(t)}\,dt+\frac 1{\sqrt{n}}dW_t,\quad t\in [0,1].
\]
The global variants are straightforward.


\subsection{Exponential observations\label{sec:EO}}

Let $X=(0,\infty )$ and $\mu (dx)$ be the Lebesgue measure on $X.$ Let $%
{\cal P}$ be the family of exponential distributions $P_\theta (dx),$ on $X$
with parameter $\theta \in \Theta =(-\infty ,0)$
\[
P_\theta (x)=e^{\theta x-V(\theta )}\mu (dx),
\]
where $V(\theta )=\log \theta .$ The corresponding Fisher information is $%
I(\theta )=\theta ^{-2},$ so condition (\ref{BKG.Ass 1}) holds true with $%
\Theta _0=[\theta _{\min },\theta _{\max }],$ for some constants $-\infty
<\theta _{\min }<\theta _{\max }<0$ and $\varepsilon _0$ small enough. Then
the parameter set $\Sigma $ contains all the functions $f(t)$ from the
H\"{o}lder ball ${\cal H}(\beta ,L),$ which satisfy $\theta _{\min }\leq
f(t)\leq \theta _{\max }.$

Consider a sequence of independent exponential observations
\begin{equation}
X_i=X(\theta _i),\quad i=1,...,n,  \label{Ex.7}
\end{equation}
with canonical parameter $\theta _i=f(t_i),$ $t_i=i/n,$ where unknown
function $f$ is in the set $\Sigma .$ Let $f_0\in \Sigma $ and ${\cal E}%
_{f_0}^n$ be the local experiment generated by the observations (\ref{Ex.7}%
). Then, by Theorem \ref{T L1}, the experiment ${\cal E}_{f_0}^n$ is
asymptotically equivalent to the local Gaussian experiment ${\cal G}_{f_0}^n$
generated by the observations
\[
dY_t^n=f(t)dt+\frac 1{\sqrt{n}}f_0(t)dW_t,\quad t\in [0,1],\quad f\in \Sigma
_{f_0}(\gamma _n).
\]

A variance-stable form for the global experiment can be also obtained. For
this we remark that $b(\theta )=-\theta $ and $F(\lambda )=\log \lambda .$
Thus, by Theorem \ref{T VS1}, a variance-stable form (local or global) is
given by the equation
\[
dY_t^n=\log \left( -f(t)\right) dt+\frac 1{\sqrt{n}}dW_t,\quad t\in [0,1].
\]


\subsection{Application to the density model\label{sec:DenMod}}

Assume that we observe a sequence of i.i.d.\ r.v.'s $X_i,$ $i=1,...,n,$ each
with density $f\in \Sigma ,$ where the set $\Sigma $ is as in the previous
section. By a poissonization technique, one can show that this experiment is
asymptotically equivalent to observing a sequence of Poisson r.v.'s $X_i,$ $%
i=1,...,n,$ with intensities $f(t_i),$ $i=1,...,n,$ respectively, where $%
t_i=\frac in,$ $i=1,...,n,$ is the uniform grid on the unit interval $[0,1].$
We skip this technical step, since its proof is similar to that given in
Nussbaum \cite{Nuss}, Section 4.

The conclusion we draw from this fact and from the example in Section \ref
{sec:PO} is that estimating a density $f(t)$ from i.i.d.\ data is
asymptotically equivalent to estimating $2$ times the square root of $f(t)$
in white noise. Thus we recover the main result of \cite{Nuss}.


\section{Local approximation\label{sec:PLE}}

\setcounter{equation}{0}


\subsection{Bounds for the Hellinger distance\label{sec:NotRes}}

Let $f_0\in \Sigma $ and let $\gamma _n$ be the nonparametric shrinking rate
defined by (\ref{F.8}). Recall briefly the setting from Section \ref{sec:NFP}%
. Set $t_i=\frac in,$ $i=1,...,n, $ and $P^{n}_{f_0,f}=P^{n}_{f},$ for $f\in
\Sigma_{f_0}(\gamma_n).$ Consider the local experiment
\begin{equation}
{{\cal {E}}}_{f_0}^n=\left( X^n,{\cal B}(X^n), \left\{ P^{n}_{f_0,f}:f\in
\Sigma _{f_0}(\gamma _n)\right\} \right) ,  \label{R.0}
\end{equation}
generated by the discrete independent observations
\begin{equation}
X_i=X(f(t_i)),\quad i=1,...,n,\quad f\in \Sigma _{f_0}(\gamma _n),
\label{R.1}
\end{equation}
with distributions in the exponential family ${\cal P}.$ Its accompanying
local Gaussian experiment
\begin{equation}
{\cal Y}_{f_0}^{n}=\left( R^n,{\cal B}(R^n),\left\{ G_{f_0,f}^{n}:f\in
\Sigma _{f_0}(\gamma _n)\right\} \right) ,  \label{R.1a}
\end{equation}
is generated by the independent Gaussian observations in discrete time
\begin{equation}
Y_i=f(t_i)+I\left( f_0(t_i)\right) ^{-\frac 12}\varepsilon_i, \quad
i=1,...,n, \quad f\in \Sigma _{f_0}(\gamma _n),  \label{R.2}
\end{equation}
with i.i.d.\ standard normal r.v.'s $\varepsilon _i.$

\begin{theorem}
\label{T LL1}Assume that $\beta >\frac 12.$ Then, for any $f_0\in \Sigma $
and any $n=1,2,...,$ the experiments ${\cal {E}}_{f_0}^n$ and ${\cal Y}%
_{f_0}^n$ can be constructed on the measurable space $(R^n,{\cal B}(R^n))$
such that
\[
\sup_{f_0\in \Sigma }\sup_{f\in \Sigma _{f_0}(\gamma
_n)}H^2(P_{f_0,f}^n,G_{f_0,f}^n)=O\left( n^{-\frac{2\beta -1}{2\beta +1}%
}(\log n)^{\frac{14\beta +5}{2\beta +1}}\right) .
\]
\end{theorem}

For the proof we make use of the following assertion which is stronger than
Theorem \ref{T LD1d} and which is also of independent interest. This theorem
corresponds to the local experiments obtained by looking only at
observations from a shrinking time interval, which after  rescaling leads to
neighborhoods of the ''almost $n^{-1/2}$'' size. Since it is the basic
result in the whole theory, we describe the setting in more detail.

Let $f_0\in {\cal L}$ (see Section \ref{sec:NFP}) and let $\gamma
_n^{*}=\kappa _0^{*}(n/\log n)^{-1/2}$ be the shrinking rate defined by (\ref
{PN.2}), with constant $\kappa _0^{*}$ being arbitrary positive. In the
sequel we always assume that the function $g\in \Sigma $ is so that $%
f=f_0+\gamma _n^{*}g\in {\cal L},$ which ensures that the corresponding
measures are defined. Consider the local experiment ${\cal E}_{f_0}^{*,n}$
defined in Section \ref{sec:almostpar}
\begin{equation}
{\cal E}_{f_0}^{*,n}=\left( X^n,{\cal B}(X^n),\left\{ P_{f_0,g}^{*,n}:g\in
\Sigma \right\} \right) ,  \label{R.10}
\end{equation}
which is generated by the discrete time independent observations
\begin{equation}
X_i=X(f(t_i)),\quad i=1,...,n,\quad f=f_0+\gamma _n^{*}g,\;g\in \Sigma ,
\label{R.11}
\end{equation}
with distributions in the exponential family ${\cal P}$ [cf. (\ref{PN.2a})].
Consider also the accompanying local Gaussian experiment from Section \ref
{sec:almostpar}, relation (\ref{PN.2b}),
\begin{equation}
{\cal Y}_{f_0}^{*,n}=\left( R^n,{\cal B}(R^n),\left\{ G_{f_0,g}^{*,n}:g\in
\Sigma \right\} \right) ,  \label{R.11a}
\end{equation}
which is generated by the independent Gaussian observations in discrete time
\begin{equation}
Y_i^{*}=f(t_i)+I\left( f_0(t_i)\right) ^{-\frac 12}\varepsilon _i,\quad
i=1,...,n,\quad f=f_0+\gamma _n^{*}g,\quad g\in \Sigma ,  \label{R.12}
\end{equation}
with i.i.d.\ standard normal r.v.'s $\varepsilon _i.$

\begin{theorem}
\label{T LL2}Let $\beta >\frac 12.$ Then, for any $f_0\in {\cal L}$ and any
$n=1,2,...,$ the local experiments ${\cal {E}}_{f_0}^{*,n}$ and
${{\cal {Y}}}_{f_0}^{*,n}$ can be constructed on the same measurable space
$(R^n{\bf ,}{\cal B}(R^n))$ such that
\[
\sup_{f_{0}\in {\cal L}}\sup_{g\in \Sigma
}H^2(P_{f_0,g}^{*,n},G_{f_0,g}^{*,n})=O\left( n^{-1}(\log n)^7\right) .
\]
\end{theorem}

\begin{remark}
Theorems {\rm \ref{T L1d}} and {\rm \ref{T LD1d}} follow immediately from
Theorems {\rm \ref{T LL1}} and {\rm \ref{T LL2}} respectively, by ({\rm \ref
{SE.6}}). The continuous versions of these results, presented in Theorems
{\rm \ref{T L1} }and {\rm \ref{T LD1}} can be established in the same way as
in Brown and Low \cite{Br-Low} and, therefore, will be not detailed here.
\end{remark}

The construction for Theorem \ref{T LL2} heavily relies upon the results on
strong approximation in Section \ref{sec:KMT} and is given in the next
Section \ref{sec:Con}. The proof of Theorem \ref{T LL2} is presented in
Section \ref{sec:ProofLE}.

Theorem \ref{T LL1} is a consequence of Theorem \ref{T LL2}. Its proof is
presented in Section \ref{sec:ParamLE}.


\subsection{Construction of a local experiment\label{sec:Con}}

In this section we give the construction of the experiments in Theorem \ref
{T LL2}. The construction on the same probability space of a Gaussian
experiment and of an exponential one is particularly simple in view of the
convenient form of the likelihood ratio of the latter [cf. (\ref{F.13})].
This form allows to employ the strong approximation result of Section \ref
{sec:KMT}. We proceed to describe formally this construction.

Let $f_0$ be a fixed function in the parameter set ${\cal L}.$ Consider the
local experiments ${\cal {E}}_{f_0}^{*,n}$ and ${{\cal {Y}}}_{f_0}^{*,n}$
from Theorem \ref{T LL2}. Recall that the shrinking rate $\gamma _n^{*}$ is
defined by
\[
\gamma _n^{*}=\kappa _0^{*}(n/\log n)^{-\frac 12},
\]
where $\kappa _0^{*}$ is arbitrary positive. Since $f_0$ is considered
fixed, for the sake of brevity, we will write $\widetilde{P}%
_g^n=P_{f_0,g}^{*,n}$ and $\widetilde{G}_g^n=G_{f_0,g}^{*,n}$. The
corresponding expectations are denoted by ${\bf E}_{\widetilde{P}_g^n}$ and $%
{\bf E}_{\widetilde{G}_g^n}$ respectively. Given a function $f=f_0+\gamma
_n^{*}g,$ $g\in \Sigma ,$ we consider the likelihood ratio of the local
experiment ${\cal {E}}_{f_0}^{*,n},$ which according to (\ref{F.13}) is
\[
\frac{d\widetilde{P}_g^n}{d\widetilde{P}_0^n}=\exp \left\{ \gamma
_n^{*}\sum_{i=1}^ng(t_i)U(X_i)-\sum_{i=1}^n\left(
V(f(t_i))-V(f_0(t_i))\right) \right\} ,
\]
where $U(X_i),$ $i=1,...,n,$ are independent r.v.'s. Denote by $\overline{U}%
(X_i)$ the r.v.\ $U(X_i)$ centered, under the measure $\widetilde{P}_0^n,$
i.e.\
\[
\overline{U}(X_i)=U(X_i)-V^{\prime }(f_0),
\]
so that, by (\ref{BKG.3a}) and (\ref{BKG.3b}), we have
\[
{\bf E}_{\widetilde{P}_0^n}\overline{U}(X_i)=0,\quad {\bf E}_{\widetilde{P}%
_0^n}\overline{U}(X_i)^2=I(f_0(t_i))=V^{\prime \prime }(f_0(t_i)).
\]
The corresponding local Gaussian experiment ${\cal Y}_{f_0}^{*,n}$ has
likelihood ratio
\begin{equation}
\frac{d\widetilde{G}_g^n}{d\widetilde{G}_0^n}=\exp \left\{ \gamma
_n^{*}\sum_{i=1}^ng(t_i)I(f(t_i))^{1/2}\varepsilon _i-\frac 12(\gamma
_n^{*})^2\sum_{i=1}^ng(t_i)^2I(f_0(t_i))\right\} ,  \label{C.4}
\end{equation}
with $f=f_0+\gamma _n^{*}g\in {\cal L}$ and $g\in \Sigma .$
Set ${\bf P}=\widetilde{G}_0^n$
and consider the probability space $(R^n,{\cal B}(R^n),{\bf P}),$ on which
$\varepsilon _i,$ $i=1,...,n,$ is a sequence of i.i.d.\ standard normal
r.v.'s. Set $N_i=I(f(t_i))^{1/2}\varepsilon _i.$ Thus we are given a
sequence of independent normal r.v.'s $N_i,$ $i=1,...,n,$ with zero means
and variances ${\bf E}N_i^2=I(f_0(t_i)),$ for $i=1,...,n.$ Because of
condition (\ref{BKG.Ass 1}) and Proposition \ref{P BKG1} we can apply
Theorem \ref{KMT.main-theorem}, according to which on this probability space
there is a sequence of independent r.v.'s $\widetilde{U}_i$ $i=1,...,n,$
such that $\widetilde{U}\stackrel{d}{=}\overline{U}(X_i),$ for any $%
i=1,...,n,$ and  for any function $g\in {\cal H}(\frac 12,L),$
\begin{equation}
{\bf P}\left( |S_n(g)|\geq x(\log n)^2\right) \leq c_1\exp \left\{
-c_2x\right\} ,\quad x\geq 0,  \label{C.5}
\end{equation}
where
\[
S_n(g)=\sum_{i=1}^ng(t_i)\left( \widetilde{U}_i-N_i\right) ,
\]
and $c_1,$ $c_2$ are constants depending only on $I_{\min },$ $I_{\max },$ $%
\varepsilon _0,$ $L.$

Now we proceed to construct a version of the likelihood process for ${\cal {E%
}}_{f_0}^{*,n}$ on the probability space $(R^n,{\cal B}(R^n),{\bf P}).$ For
this define the experiment ${\cal F}_{f_0}^{*,n}$ as follows:
\[
{\cal F}_{f_0}^{*,n}=\left( R^n,{\cal B}(R^n),\left\{ F_g^n:\;g\in \Sigma
\right\} \right) ,
\]
where $F_g^n$ is the probability measure on the measurable space $(R^n,{\cal %
B}(R^n))$ defined, for any $g\in \Sigma, $ by the equality
\begin{equation}
\frac{dF_g^n}{d{\bf P}}=\exp \left\{ \gamma _n^{*}\sum_{i=1}^ng(t_i)%
\widetilde{U}_i-\sum_{i=1}^n\left( V(f(t_i))-V(f_0(t_i))-\gamma
_n^{*}g(t_i)V^{\prime }(f_0(t_i))\right) \right\},  \label{C.6}
\end{equation}
with $f=f_0+\gamma _n^{*}g$ ; then $F_0^n={\bf P}.$

Let us remark that, since the sequences $\widetilde{U}_i,$ $i=1,...,n,$ and $%
\overline{U}(X_i),$ $i=1,...,n,$ have the same joint distributions and the
sufficient statistic in the experiment ${\cal {E}}_{f_0}^{*,n}$ is
\[
\sum_{i=1}^ng(t_i)U(X_i),
\]
the experiments ${\cal E}_{f_0}^{*,n}$ and ${\cal {F}}_{f_0}^{*,n}$ are
exactly equivalent. Therefore, we can assume in the sequel that the
construction is performed on the initial probability space, i.e. that $%
\widetilde{U}_i\equiv \overline{U}(X_i),$ for $i=1,...,n,$ and ${\cal E}%
_{f_0}^{*,n}\equiv {\cal {F}}_{f_0}^{*,n}.$ In particular $\widetilde{P}%
_g^n\equiv F_g^n,$ for any $g\in \Sigma ,$ and $\widetilde{P}_0^n\equiv {\bf %
P,}\equiv \widetilde{G}_0^n.$


\subsection{Proof of local equivalence: almost $n^{-1/2}$-neighborhoods\label%
{sec:ProofLE}}

In this section we present a proof of Theorem \ref{T LL2}. First recall that
according to the last remark in the previous section we consider the
experiment ${\cal {E}}_{f_0}^{*,n}$ to be a version of the original one
constructed such that its likelihood process is on a common probability
space $(R^n,{\cal B}(R^n),{\bf P})$ with the Gaussian likelihood process for
${\cal Y}_{f_0}^{*,n}.$ Recall also that the "central" measures (i.e.\ those
with $g=0$) in the local experiments ${\cal E}_{f_0}^{*,n}$ and ${\cal Y}%
_{f_0}^{*,n}$ coincide with the measure ${\bf P}$. We will denote the
expectation under ${\bf P}$ by ${\bf E .}$

Let $f_0\in {\cal L}$ and $f=f_0+\gamma _n^{*}g,$ $g\in \Sigma .$ Note that
the function $g$ belongs also to the H\"{o}lder ball ${\cal H}(\frac 12,L),$
since $\beta \geq \frac 12.$ Then, by taking $x=\frac{c_3}{c_2}\log n$ in
(\ref{C.5}), we arrive at
\begin{equation}
{\bf P}(|S_n(g)|\geq \frac{c_3}{c_2}\gamma _n(\log n)^3)\leq c_1\exp
\{-c_3\log n\},  \label{DLE.2}
\end{equation}
where $c_1,$ $c_2$ are the same as in (\ref{C.5}) and $c_3$ is a ''free''
constant whose value will be chosen later. Recall that according to our
agreement $\widetilde{U}_i=\overline{U}(X_i),$ $i=1,...,n,$ and therefore
\[
S_n(g)=\sum_{i=1}^ng(t_i)(\overline{U}(X_i)-N_i).
\]
What we have to prove is that the Hellinger distance between the measures $%
\widetilde{P}_g^n$ and $\widetilde{G}_g^n$ satisfies
\begin{equation}
H^2\left( \widetilde{P}_g^n,\widetilde{G}_g^n\right) \leq c_4n^{-1}(\log
n)^7,  \label{DL.2}
\end{equation}
for some constant $c_4$ depending only on $I_{\max },$ $I_{\min },$ $%
\varepsilon _0,$ $\kappa _0^{*},$ $L.$

Well known properties of the Hellinger distance [see (\ref{SE.1}] and (\ref
{SE.6a})] imply
\[
H^2\left( \widetilde{P}_g^n,\widetilde{G}_g^n\right) =H^2\left( \Lambda
^1(g),\Lambda ^2(g)\right) =\frac 12{\bf E}\left( \Lambda ^1(g)^{\frac
12}-\Lambda ^2(g)^{\frac 12}\right) ^2,
\]
where we denote for brevity
\[
\Lambda ^1(g)=\frac{d\widetilde{P}_g^n}{d{\bf P}},\quad \Lambda ^2(g)=\frac{d%
\widetilde{G}_g^n}{d{\bf P}}.
\]
Set $u_n=\frac{c_3}{c_2}\gamma _n(\log n)^3.$ Then obviously
\[
H^2\left( \widetilde{P}_g^n,\widetilde{G}_g^n\right) \leq J_1+J_2,
\]
where [${\bf 1}(A)$ being the indicator for event $A$]
\begin{eqnarray*}
J_1 &=&\frac 12{\bf E}\,{\bf 1}\left( |S_n(g)|<u_n\right) \left( \Lambda
^1(g)^{\frac 12}-\Lambda ^2(g)^{\frac 12}\right) ^2, \\
J_2 &=&\frac 12{\bf E}\,{\bf 1}\left( |S_n(g)|\geq u_n\right) \left( \Lambda
^1(g)^{\frac 12}-\Lambda ^2(g)^{\frac 12}\right) ^2.
\end{eqnarray*}
First we give an estimate for $J_1$. Changing the probability measure we
obtain
\[
J_1=\frac 12{\bf E}_{\widetilde{G}_g^n}{\bf 1}\left( |S_n(g)|\leq u_n\right)
\left( \Lambda ^2(g)^{\frac 12}\Lambda ^1(g)^{-\frac 12}-1\right) ^2.
\]
According to (\ref{C.6}) and (\ref{C.4})
\begin{equation}
\Lambda ^1(g)=\exp \left\{ \gamma
_n^{*}\sum_{i=1}^ng(t_i)U(X_i)-\sum_{i=1}^n\left(
V(f(t_i))-V(f_0(t_i))\right) \right\}  \label{DLE.2a}
\end{equation}
and
\[
\Lambda ^2(g)=\exp \left\{ \gamma _n^{*}\sum_{i=1}^ng(t_i)N_i-\frac
12(\gamma _n^{*})^2\sum_{i=1}^ng(t_i)^2I(f_0(t_i))\right\} .
\]
This gives
\[
\frac{\Lambda ^2(g)}{\Lambda ^1(g)}=\exp \left\{ -S_n(f)+R(f_0,f)\right\} ,
\]
where
\begin{eqnarray*}
R(f_0,f) &=&\sum_{i=1}^n\{V(f(t_i))-V(f_0(t_i))-\gamma _n^{*}g(t_i)V^{\prime
}(f_0(t_i)) \\
&&\quad \quad \quad \quad \;\;\;\;\;\;-\frac 12(\gamma
_n^{*})^2g(t_i)^2V^{\prime \prime }(f_0(t_i))\}.
\end{eqnarray*}
A three term Taylor expansion yields
\[
|R(f_0,f)|\leq \frac 16(\gamma _n^{*})^3\sum_{i=1}^n\left| g(t_i)\right|
^3\left| V^{\prime \prime \prime }(\widetilde{f}_i)\right| ,
\]
with $\widetilde{f}_i=f_0(t_i)+\nu _i\gamma _n^{*}g(t_i),$ $0\leq \nu _i\leq
1.$ Since $||g||_\infty \leq L$ and $\left| V^{\prime \prime \prime }(%
\widetilde{f}_i)\right| \leq c_5,$ with $c_5$ depending only on $I_{\max },$
$\varepsilon _0$ (see Proposition \ref{P BKG2}), it is clear that
\[
|R(f_0,f)|\leq \frac 16L^3c_5(\gamma _n^{*})^3n\leq c_6n^{-\frac 12}(\log
n)^{\frac 32},
\]
with $c_6$ depending on $c_5,$ $L,$ $\kappa _0^{*}$ and for $n$ large
enough. This yields the following estimate:
\begin{eqnarray*}
J_1 &=&\frac 12{\bf E}_{\widetilde{G}_f^n}{\bf 1}(|S_n(f)|\leq u_n)\left(
\exp \left\{ -\frac 12S_n(f)+\frac 12R(f_0,f)\right\} -1\right) ^2 \\
\ &\leq &c_7\left( \left| S_n(f)\right| +\left| R(f_0,f)\right| \right) ^2 \\
\ &\leq &c_7\left( \frac{c_2}{c_1}n^{-\frac 12}(\log n)^{\frac
72}+c_6n^{-\frac 12}(\log n)^{\frac 32}\right) ^2,
\end{eqnarray*}
with $c_7$ depending on $c_2,$ $c_3,$ $c_6.$ Hence with some $c_8$ depending
on $c_2,$ $c_3,$ $c_6,$ $c_7$
\begin{equation}
J_1\leq c_7n^{-1}(\log n)^7.  \label{DL.7}
\end{equation}

Now we proceed to estimate $J_2.$ The H\"{o}lder inequality implies
\begin{equation}
J_2\leq \frac 12J_3^{\frac 12}J_4^{\frac 12},  \label{DLE.3}
\end{equation}
where by (\ref{DLE.2})
\begin{equation}
J_3={\bf P}(|S_n(f)|\geq u_n)\leq c_1\exp \{-c_3\log n\},  \label{DL.8}
\end{equation}
for $n$ large enough, and
\begin{equation}
J_4={\bf E}\left( \Lambda ^1(g)^{\frac 12}+\Lambda ^2(g)^{\frac 12}\right)
^4.  \label{DL.8a}
\end{equation}
Note that the constant $c_3$ in (\ref{DL.8}) is ''free''. We will show that $%
J_4$ is bounded from above by $32\exp \{c_9\log n\}$ with some constant $%
c_9$ which we calculate below. Indeed, from (\ref{DL.8a}) we get
\begin{equation}
J_4\leq 16\left( {\bf E}(\Lambda ^1(g))^2+{\bf E}(\Lambda ^2(g))^2\right) .
\label{DL.8b}
\end{equation}
First we give a bound for ${\bf E}(\Lambda ^1(g))^2.$ It follows from (\ref
{DLE.2a}) that
\[
{\bf E}(\Lambda ^1(g))^2={\bf E}\exp \{2\gamma _n^{*}\sum_{i=1}^ng(t_i)%
\overline{U}(X_i)-2R_0(f_0,f)\},
\]
where
\[
R_0(f_0,f)=\sum_{i=1}^n\left\{ V(f(t_i))-V(f_0(t_i))-\gamma
_n^{*}g(t_i)V^{\prime }(f_0(t_i))\right\} .
\]
The estimate for the remainder $R_0(f_0,f)$ is straightforward by Taylor's
formula:
\[
|R_0(f_0,f)|\leq \frac 12(\gamma _n^{*})^2\sum_{i=1}^n\left| g(t_i)\right|
^2\left| V^{\prime \prime }(\widetilde{f}_i)\right| \leq \frac 12L^2I_{\max
}(\gamma _n^{*})^2n,
\]
where $\widetilde{f}_i=f(t_i)+\nu _i\gamma _n^{*}g(t_i),$ $0\leq \nu _i\leq
1.$ Since $\overline{U}(X_i),$ $i=1,...,n,$ are independent r.v.'s, using
Proposition \ref{P BKG1} we obtain
\begin{eqnarray*}
{\bf E}(\Lambda ^1(g))^2 &=&\exp \left\{ -2R_0(f_0,f)\right\} \prod_{i=1}^n%
{\bf E}\exp \left\{ 2\gamma _n^{*}g(t_i)\overline{U}(X_i)\right\} \\
\ &\leq &\exp \left\{ 2|R_0(f_0,f)|\right\} \prod_{i=1}^n\exp \left\{
2(\gamma _n^{*})^2g(t_i)^2I_{\max }\right\} \\
\ &\leq &\exp \left\{ 3L^2I_{\max }(\gamma _n^{*})^2n\right\} \leq \exp
\left\{ c_9\log n\right\} ,
\end{eqnarray*}
where $c_9$ depends on $I_{\max },$ $L,$ $\kappa _0^{*}.$ The bound ${\bf E}%
(\Lambda ^1(g))^2\leq \exp \{c_9\log n\}$ can be proved similarly. These
bounds and (\ref{DL.8b}) yield
\begin{equation}
J_4\leq 32\exp \left\{ c_9\log n\right\} .  \label{DL.9}
\end{equation}
Using the bounds for $J_3$ and $J_4$ given by (\ref{DL.8}) and (\ref{DL.9})
in (\ref{DLE.3}) we obtain
\[
J_2\leq 2\sqrt{2c_1}\exp \left\{ -\frac 12(c_3-c_9)\log n\right\} ,
\]
from which, taking $c_3=c_9+2,$ we get
\begin{equation}
J_2\leq 2\sqrt{2c_1}/n\leq 2\sqrt{2c_1}n^{-1}(\log n)^7.  \label{DL.10}
\end{equation}
The desired inequality (\ref{DL.2}) follows from (\ref{DL.7}) and (\ref
{DL.10}).


\subsection{Proof of local equivalence: nonparametric neighborhoods\label%
{sec:ParamLE}}

We present here a proof for Theorem \ref{T LL1}. Before the rigorous
argument let us briefly expound the main idea. We start by splitting the
original local experiment ${{\cal {E}}}_{f_0}^n$ into $m$ parts that
correspond to fractions of the observations over shrinking time intervals
having length of order $\delta _n=\gamma _n^{1/\beta },$ where $\gamma _n$
is the shrinking rate of the neighborhood $\Sigma _{f_0}(\gamma _n).$ One
may call the corresponding experiments {\it doubly local}. Denote by $%
n_k=O(n\delta _n)$ the number of observations (i.e.\ number of design points
$t_i$) in the $k$-th doubly local experiment. After rescaling this one can
be viewed as an experiment on the whole interval $[0,1]$ over a shrinking
neighborhood of size $O(n_k/\log n_k)^{-\frac 12}.$ By Theorem \ref{T LL2}
we can ''approximate'' this experiment by the corresponding Gaussian
experiment, with a bound for the {\it squared} Hellinger distance between
corresponding measures of order $O(n_k^{-1}(\log n_k)^7).$ Further arguments
are based on the crucial inequality (\ref{SE.2}), which is applied to the
original (parameter-local) experiment on the unit interval $[0,1]$
construed as a product of the $m$ local doubly local experiments on the
intervals of size $\delta _n$. Since the Gaussian experiment ${\cal Y}%
_{f_0}^n$ can be decomposed similarly, we obtain a bound for the squared
Hellinger distance between ${\cal E}_{f_0}^n$ and ${\cal Y}_{f_0}^n$
\[
O(mn_k^{-1}(\log n_k)^7)=O(n\delta _n^{-2}(\log n_k)^7)=o(1)
\]
as $n\rightarrow \infty ,$ for $\beta >\frac 12,$ which proves our theorem.

Now we turn to the argument in detail. Let $\beta >\frac 12$ and $f_0\in
\Sigma .$ Assume that the shrinking rate of the neighborhood $\Sigma
_{f_0}(\gamma _n)$ is given by
\begin{equation}
\gamma _n=\kappa _0(n/\log n)^{-\frac \beta {2\beta +1}}  \label{P.1}
\end{equation}
with some constant $\kappa _0$ depending on $\beta .$ Set $m=[1/\delta _n],$
where
\begin{equation}
\delta _n=(\gamma _n)^{\frac 1\beta }=\kappa _0^{1/\beta }(n/\log n)^{-\frac
1{2\beta +1}}.  \label{ParLE.1}
\end{equation}
Consider a partition ${\frak A}$ of the unit interval $[0,1]$ into intervals
$A_k=(a_k,b_k],$ $k=1,...,m,$ of length $1/m.$ It is easy to see that for $n$
large enough
\begin{equation}
\frac 1{2\delta _n}\leq m\leq \frac 1{\delta _n}.  \label{ParLE.2}
\end{equation}
Set $I_k=\{i:t_i\in A_k\}.$ Denote by $n_k$ the cardinality of $I_k,$ i.e.\
the number of design points which fall into the interval $A_k,$ $k=1,...,m.$
It is clear that for $n$ large enough
\begin{equation}
n\delta _n\leq n_k\leq 2n\delta _n.  \label{ParLE.3}
\end{equation}
We particularly point out that (\ref{P.1}), (\ref{ParLE.1}), (\ref{ParLE.2})
and (\ref{ParLE.3}) imply
\begin{equation}
\gamma _n\leq \gamma _n^{*}\equiv \kappa _0^{*}(n_k/\log n_k)^{-1/2}
\label{P.4}
\end{equation}
for some constant $\kappa _0^{*}$ depending on $\kappa _0$ and $\beta .$ Let
$a_k(t)$ be the linear function that maps the unit interval $[0,1]$ onto the
interval $[a_k,b_k],$ i.e.\ $a_k(t)=t/m+a_k,$ $t\in [0,1].$ For any $f\in
\Sigma _{f_0}(\gamma _n)$ and $k=1,...,m,$ consider the function $f_k$
defined on the interval $[0,1]$ as follows:
\begin{equation}
f_k(t)=\frac{(f-f_0)(a_k(t))}{\gamma _n^{*}}.  \label{P.4a}
\end{equation}
We will prove that $f_k\in \Sigma ^{*}=\Sigma (\beta ,L^{*}),$ with some $%
L^{*}$ depending on $L$ and $\beta .$ Indeed, since $f$ $\in \Sigma
_{f_0}(\gamma _n),$ we have  $||f_k||_\infty \leq ||f-f_0||_\infty /\gamma
_n^{*}\leq \gamma _n/\gamma _n^{*}\leq 1.$ On the other hand, since $f,f_0$
are in the H\"{o}lder ball ${\cal H}(\beta ,L),$ the function $\psi =f-f_0$
is also in the H\"{o}lder ball. Taking into account (\ref{ParLE.1}), (\ref
{ParLE.2}) and (\ref{P.4}) we obtain for any $x,y\in [0,1],$%
\begin{eqnarray*}
\left| f_k^{([\beta ])}(x)-f_k{}^{([\beta ])}(y)\right| &=&m^{-[\beta
]}\left| \psi ^{([\beta ])}(a_k(x))-\psi ^{([\beta ])}(a_k(y))\right|
/\gamma _n^{*} \\
\ &\leq &m^{-[\beta ]}L\left| a_k(x)-a_k(y)\right| ^{\beta -[\beta ]}/\gamma
_n^{*} \\
\ &\leq &2^\beta \delta _n^\beta L\left| x-y\right| ^{\beta -[\beta
]}/\gamma _n^{*} \\
\ &=&2^\beta \gamma _nL\left| x-y\right| ^{\beta -[\beta ]}/\gamma _n^{*} \\
\ &\leq &2^\beta L\left| x-y\right| ^{\beta -[\beta ]},
\end{eqnarray*}
proving that $f_k$ is in the H\"{o}lder ball $\Sigma (\beta ,L^{*})$ with $%
L^{*}=2^\beta L.$

Let $X^{k,n}=\left\{ X_i,\;i\in I_k\right\} $ be the fragment of
observations $\left\{ X_i,\;i=1,...,n\right\} $ [defined by (\ref{R.1})]
associated to the time interval $A_k,$ for some $k\in \{1,...,m\}.$ After a
rescaling with the linear function $a_k(t)$ these observations can be
associated to design points $t_i^k=\frac i{n_k},$ $i=1,...,n_k,$ on the unit
interval $[0,1].$ Let $P_f^{n,k}$ be the measure on $(R^{n_k},{\cal B}%
(R^{n_k}))$ induced by the set of r.v.'s $X^{n,k}$ and set $%
P_{f_0,g}^{k,n}=P_f^{k,n}$ for $f=f_0+\gamma _n^{*}g$. For each $k\in
\{1,...,m\}$ consider the local experiment
\begin{equation}
{{\cal {E}}}_{f_0}^{k,n}=\left( R^{n_k},{\cal B}(R^{n_k}),\left\{
P_{f_0,g}^{k,n}:\;g\in \Sigma ^{*}\right\} \right) .   \label{P.4b}
\end{equation}
In the same way we introduce the local Gaussian experiment
\begin{equation}
{{\cal {Y}}}_{f_0}^{k,n}=\left( R^{n_k},{\cal B}(R^{n_k}),\left\{
G_{f_0,g}^{k,n}:\;g\in \Sigma ^{*}\right\} \right)   \label{P.4c}
\end{equation}
generated by the observation fragment $Y^{k,n}=\left\{ Y_i,\;i\in
I_k\right\} ,$ with $Y_i$ defined by (\ref{L.3}). Here $G_{f_0,g}^{k,n}$ is
the Gaussian shift measure on $(R^{n_k},{\cal B}(R^{n_k}))$ induced by the
observations $Y^{n,k}$ under $f=f_0+\gamma _n^{*}g$.

According to Theorem \ref{T LL2} the experiments ${{\cal {E}}}_{f_0}^{k,n}$
and ${{\cal {Y}}}_{f_0}^{k,n}$ can be constructed on the measurable space $%
(R^{n_k},{\cal B}(R^{n_k}))$ such that
\begin{equation}
\sup_{g\in \Sigma ^{*}}H^2\left( P_{f_0,g}^{k,n},G_{f_0,g}^{k,n}\right) \le
c_1n_k^{-1}(\log n_k)^7,  \label{P.5}
\end{equation}
with a constant $c_1$ depending on $I_{\max },$ $I_{\min },$ $\varepsilon
_0, $ $\kappa _0^{*},$ $L.$ Now consider the subexperiments of ${{\cal {E}}}%
_{f_0}^{k,n}$ and ${{\cal {Y}}}_{f_0}^{k,n}$ obtained by setting $g=f_k,$ $%
f\in \Sigma _{f_0}(\gamma _n)$ in (\ref{P.4b}) and (\ref{P.4c}), where $f_k$
is defined by (\ref{P.4a}). Reindex those subexperiments by $f\in \Sigma
_{f_0}(\gamma _n)$ and call them $\widetilde{{\cal E}}_{f_0}^{k,n},$ $%
\widetilde{{\cal Y}}_{f_0}^{k,n}$ respectively. For the respective reindexed
measures $\widetilde{P}_{f_0,f}^{k,n}=P_{f_0,f_k}^{k,n}$ and $\widetilde{G}%
_{f_0,f}^{k,n}=G_{f_0,f_k}^{k,n}$ we have as a consequence of (\ref{P.5})
\begin{equation}
\sup_{f\in \Sigma _{f_0}(\gamma _n)}H^2\left( \widetilde{P}_{f_0,f}^{k,n},%
\widetilde{G}_{f_0,f}^{k,n}\right) \le c_1n_k^{-1}(\log n_k)^7,  \label{P.5a}
\end{equation}
since $f_k\in \Sigma ^{*},$ for any $f\in \Sigma _{f_0}(\gamma _n).$

Define the experiment ${\Bbb E}_{f_0}^n$ as
\[
{\Bbb E}_{f_0}^n=\widetilde{{\cal E}}_{f_0}^{1,n}\otimes ...\otimes
\widetilde{{\cal E}}_{f_0}^{m,n};
\]
this one is  obviously  (exactly) equivalent to ${\cal E}_{f_0}^{*,n}$ defined by (%
\ref{R.10}), (\ref{R.11}). The corresponding local Gaussian experiment $%
{\cal Y}_{f_0}^{*,n}$ defined by (\ref{R.11a}), (\ref{R.12}) is (exactly)
equivalent to the experiment
\[
{\Bbb Y}_{f_0}^n=\widetilde{{\cal Y}}_{f_0}^{1,n}\otimes ...\otimes
\widetilde{{\cal Y}}_{f_0}^{m,n}.
\]
It remains only to note that $P_{f_0,f}^n=\widetilde{P}_{f_0,f}^{1,n}\otimes
...\otimes \widetilde{P}_{f_0,f}^{m,n}$ and $G_{f_0,f}^n=\widetilde{G}%
_{f_0,f}^{1,n}\otimes ...\otimes \widetilde{G}_{f_0,f}^{m,n},$ and thus
according to (\ref{SE.2})
\begin{eqnarray*}
H^2\left( P_{f_0,f}^n,G_{f_0,f}^{*,n}\right) &\leq &\sum_{i=1}^mH^2\left(
\widetilde{P}_{f_0,f}^{k,n},\widetilde{G}_{f_0,f}^{k,n}\right) \\
\ &\leq &mn_k^{-1}(\log n_k)^7 \\
\ &\leq &n^{-1}\delta _n^{-2}(\log n)^7 \\
\ &\leq &c_2n^{-\frac{2\beta -1}{2\beta +1}}(\log n)^{\frac{14\beta +5}{%
2\beta +1}},
\end{eqnarray*}
where $c_2$ is a constant depending on $I_{\max },$ $I_{\min },$ $%
\varepsilon _0,$ $\beta ,$ $L.$ This completes the proof of Theorem \ref{T
LL1}.


\subsection{Proof of the variance-stable form\label{sec:P-VS}}

In this section we present a proof of Theorem \ref{T VS1}.

Introduce the following experiments: ${\cal F}_{f_0}^{1,n}$ generated
by the observations
\begin{equation}
dY_t^{1,n}= (f(t)-f_0(t))I(f_0(t))^{1/2}dt+\frac 1{\sqrt{n}}dW_t, \quad
t\in [0,1],  \label{A.1}
\end{equation}
where $f\in \Sigma _{f_0}(\gamma _n),$ and ${\cal F}_{f_0}^{2,n}$ generated
by the observations
\begin{equation}
dY_t^{2,n}=\left( \Gamma \left( f(t)\right) -\Gamma \left( f_0(t)\right)
\right) dt+\frac 1{\sqrt{n}}dW_t,\quad t\in [0,1],  \label{A.2}
\end{equation}
where $f\in \Sigma _{f_0}(\gamma _n).$ We prove the following assertion.

\begin{proposition}
\label{P AAA1} The experiments ${\cal F}_{f_0}^{1,n}$ and ${\cal F}%
_{f_0}^{2,n}$ are asymptotically equivalent. Moreover
\[
\Delta ^2\left( {\cal F}_{f_0}^{1,n},{\cal F}_{f_0}^{2,n}\right) \leq c_0n^{-%
\frac{2\beta -1}{2\beta +1}}(\log n)^{\frac{4\beta }{2\beta +1}},
\]
where $c_0$ is a constant depending only on $I_{\max },$ $I_{\min }$ and $%
\varepsilon _0.$
\end{proposition}

{\noindent{\bf Proof.}\ \ } By Taylor expansion we have for any $\theta
,\theta _0\in \Theta _0,$
\begin{equation}
\Gamma (\theta )-\Gamma (\theta _0)=(\theta -\theta _0)\Gamma ^{\prime
}(\theta _0)+\frac 12(\theta -\theta _0)^2\Gamma ^{\prime \prime }(%
\widetilde{\theta }_0),  \label{A.6}
\end{equation}
where $\widetilde{\theta }_0=\theta _0+\nu (\theta -\theta _0),$ $0\leq \nu
\leq 1.$ Next it follows from (\ref{VS.001}) that
\begin{equation}
\Gamma ^{\prime }(\theta _0)=\sqrt{I(\theta _0)},  \label{A.7}
\end{equation}
while, using assumption (\ref{BKG.Ass 1}), the second derivative of $\Gamma
(\lambda )$ can easily be seen to satisfy
\begin{equation}
\left| \Gamma ^{\prime \prime }(\widetilde{\theta }_0)\right| \leq \frac
12I_{\min }^{-1/2}.  \label{A.8}
\end{equation}
Hence by (\ref{A.6}), (\ref{A.7}) and (\ref{A.8}) we get for $f\in \Sigma
_{f_0}(\gamma _n),$
\begin{equation}
\Gamma (f(t))-\Gamma (f_0(t))=\left( f(t)-f_0(t)\right)
I(f_0(t))^{1/2}+\frac 14\nu I_{\min }^{-1/2}\gamma _n^2,  \label{A.9}
\end{equation}
with some $\left| \nu \right| \leq 1.$ Set for brevity
\begin{eqnarray*}
m_1(t) &=&\left( f(t)-f_0(t)\right) I(f_0(t))^{1/2}, \\
m_2(t) &=&\Gamma (f(t))-\Gamma (f_0(t)).
\end{eqnarray*}
Let $Q_{f_0,f}^{1,n}$ and $Q_{f_0,f}^{2,n}$ be the measures induced by
observations (\ref{A.1}) and (\ref{A.2}) respectively. Then by formula (\ref
{SE.7}) and by (\ref{A.9}) we get
\begin{eqnarray*}
H^2\left( Q_{f_0,f}^{1,n},Q_{f_0,f}^{2,n}\right) &=&1-\exp \left\{ -\frac
n8\int_0^1\left( m_1(t)-m_2(t)\right) ^2dt\right\} \\
\ &\leq &\frac n8\gamma _n^4c_1=c_2n^{-\frac{2\beta -1}{2\beta +1}}(\log n)^{%
\frac{4\beta }{2\beta +1}}.\;{\mbox{\ \rule{.1in}{.1in}}}
\end{eqnarray*}

Theorem \ref{T VS1} can be obtained easily from Proposition \ref{P AAA1} if
we note that $\Delta ({\cal G}_{f_0}^n,{\cal F}_{f_0}^{1,n})=0$ and also $%
\Delta ({\cal F}_{f_0}^{2,n},\widehat{{\cal G}}_{f_0}^n)=0,$ i.e.\ these
experiments are (exactly) equivalent by the remark following immediately
after formula (\ref{SE.6a}).


\section{Global approximation\label{sec:GLOB}}

\setcounter{equation}{0}


\subsection{The preliminary estimator\label{sec:PRES}}

The following lemma provides the preliminary estimator that is necessary for
the globalization procedure over the parameter set $\Sigma $.

\begin{lemma}
\label{L ES1}Let $\beta \in (\frac 12,1).$ In the experiment ${\cal E}^n$
there is an estimator $\widehat{f}_n:X^n\rightarrow \Sigma $ taking a finite
number of values and fulfilling
\[
\sup_{f\in \Sigma }P_f^n\left( \left\| \widehat{f}_n-f\right\| _\infty
>c_1\gamma _n\right) \leq c_2\frac 1{\sqrt{n}},
\]
where $c_1$ and $c_2$ are constants depending only on $I_{\max },$ $I_{\min
},$ $\kappa _0,$ $L,$ $\beta .$
\end{lemma}

{\noindent{\bf Proof.}\ \ } Let $\gamma _n$ be given by (\ref{F.8}) and $%
\delta _n=\gamma _n^{1/\beta }.$ Introduce the kernel $K(u)$ as a bounded
function of finite support such that
\begin{equation}
0\leq K(u)\leq k_{\max },\quad K(u)=0,\ u\notin (-\tau ,\tau ),\quad
\int_{-\tau }^\tau K(u)du=1,  \label{ES.1}
\end{equation}
where $k_{\max }$ and $\tau $ are some absolute constants. We will assume
also that $K(u)$ satisfies a H\"{o}lder condition with exponent $\beta $.
Let
\begin{equation}
\rho _n(t)=\frac 1{n\delta _n}\sum_{i=1}^nK\left( \frac{t_i-t}{\delta _n}%
\right)   \label{ES.1a}
\end{equation}
where $t_i=i/n,$ $i=1,...,n.$ It is easy to see that there are two constants
$\rho _{\min }$ and $\rho _{\max }$ such that for $n$ large enough
\begin{equation}
0<\rho _{\min }\leq \rho _n(t)\leq \rho _{\max }<\infty   \label{ES.3}
\end{equation}
for any $t\in [0,1].$

Consider the functions $f\in \Sigma $ and $g(t)=b\left( f(t)\right) ,$ $t\in
[0,1],$ where $b(\theta )=V^{\prime }(\theta ),$ $\theta \in \Theta _0$ (see
Section \ref{sec:VS}). Define an estimator $g_n^{*}$ of $g$ as follows: for
any $t\in [0,1],$ set
\[
g_n^{*}(t)=\frac 1{n\delta _n\rho _n(t)}\sum_{i=1}^nK\left( \frac{t_i-t}{%
\delta _n}\right) U(X_i),
\]
where $U(x)$ is the sufficient statistic in the exponential experiment $%
{\cal E}.$ The estimator $g_n^{*}$ is known as the Nadaraya-Watson estimator. We
will show that there are two constants $c_3$ and $c_4$ depending only on $%
I_{\max },$ $\kappa _0,$ $L,$ $\beta ,$ $k_{\max },$ $\rho _{\min },$ $\tau
, $ such that
\begin{equation}
\sup_{f\in \Sigma }P_f^n\left( \left\| g_n^{*}-g\right\| _\infty >c_3\gamma
_n\right) \leq c_4\frac 1n.  \label{ES.5}
\end{equation}
First we note that by (\ref{ES.1a}),
\begin{equation}
E_f^ng_n^{*}(t)-g(t)=\frac 1{n\delta _n\rho _n(t)}\sum_{i\in J_n(t)}K\left(
\frac{t_i-t}{\delta _n}\right) \left( g(t_i)-g(t)\right) ,  \label{ES.6}
\end{equation}
for any $t\in [0,1],$ where $J_n(t)=\left\{ i:t_i\in (t-\tau \delta
_n,t+\tau \delta _n)\right\} $ and $\#J_n(t)\leq 2\tau n\delta _n.$ It is
easy to see that since $f\in \Sigma ,$ we have for $i\in J_n(t)$
\begin{equation}
\left| g(t_i)-g(t)\right| \leq I_{\max }L(2\tau \delta _n)^\beta =(2\tau
)^\beta I_{\max }L\gamma _n.  \label{ES.7}
\end{equation}
From (\ref{ES.1}), (\ref{ES.6}) and (\ref{ES.7}) we have
\begin{equation}
\left\| E_f^ng_n^{*}-g\right\| _\infty \leq c_5\gamma _n,  \label{ES.8}
\end{equation}
with some constant $c_5$ depending on $I_{\max },$ $L,$ $\beta ,$ $k_{\max
}, $ $\rho _{\min },$ $\tau .$ To handle the difference $%
g_n^{*}-E_f^ng_n^{*} $ we remark that
\begin{equation}
g_n^{*}(t)-E_f^ng_n^{*}(t)=\frac 1{n\delta _n\rho _n(t)}\sum_{i\in
J_n(t)}K\left( \frac{t_i-t}{\delta _n}\right) \overline{U}(X_i),
\label{ES.8a}
\end{equation}
where $\overline{U}(X_i)=U(X_i)-E_f^nU(X_i)=U(X_i)-g(t_i).$
Set for brevity
$L_i(t)=K\left( \frac{t_i-t}{\delta _n}\right) .$ Define a piecewise constant
approximation of $L_i(t)$ as follows: put $\widetilde{L}_i(t)=L_i(s_k)$ for $%
t\in A_k,$ where $A_1=[0,s_1],$ $A_k=(s_{k-1},s_k],$ $k=2,...,n^2,$ $%
s_k=k/n^2,$ $k=0,...,n^2.$ Since the function $K(u)$ satisfies a H\"{o}lder
condition with exponent $\beta ,$ there is a constant $c_6$ such that
\begin{equation}
\left\| L_i-\widetilde{L}_i\right\| _\infty \leq c_6\left( \frac 1n\right)
^{2\beta }.  \label{ES.8b}
\end{equation}
Then (\ref{ES.8a}) and (\ref{ES.3}) imply
\[
\left\| g_n^{*}-E_f^ng_n^{*}\right\| _\infty \leq Q_1+Q_2,
\]
where
\begin{eqnarray*}
Q_1 &=&\frac 1{n\delta _n\rho _{\min }}\sup_{t\in [0,1]}\left|
\sum_{i=1}^n\left( L_i(t)-\widetilde{L}_i(t)\right) \overline{U}(X_i)\right|
, \\
Q_2 &=&\frac 1{n\delta _n\rho _{\min }}\sup_{t\in [0,1]}\left| \sum_{i=1}^n%
\widetilde{L}_i(t)\overline{U}(X_i)\right| .
\end{eqnarray*}
Using (\ref{ES.8b}) we get
\[
Q_1\leq \frac{c_6}{n^{2\beta +1}\delta _n\rho _{\min }}\sum_{i=1}^n\left|
\overline{U}(X_i)\right| .
\]
Set for brevity $u_n=c_6\left( \rho _{\min }n^{2\beta +1}\delta _n\gamma
_n\right) ^{-1}\log n.$ Then, with a ''free'' constant $c_7\geq 1,$
\[
P_f^n(Q_1>c_7\gamma _n)\leq e^{-c_7\log n}\prod_{i=1}^nE_f^n\exp \left\{
u_n\left| \overline{U}(X_i)\right| \right\} .
\]
Since $u_n\leq \varepsilon _0,$ for $n$ large enough it is easy to see
using Proposition \ref{P BKG1} that
\[
E_f^n\exp \left\{ u_n\left| \overline{U}(X_i)\right| \right\} \leq \exp
\left\{ u_nc_8\right\} ,
\]
where $c_8$ is a constant depending only on $I_{\max },$ $\varepsilon _0.$
As $nu_n\rightarrow 0$ for $n\rightarrow \infty ,$ we have for
sufficiently large $n$%
\begin{equation}
P_f^n(Q_1>c_7\gamma _n)\leq \exp \left\{ -c_7\log n+c_8nu_n\right\} \leq
2\exp \left\{ -c_7\log n\right\} \leq 2\frac 1n.  \label{ES.8d}
\end{equation}
To obtain a bound for $Q_2,$ we remark that $\widetilde{L}_i(t)$ is
piecewise constant and $\widetilde{L}_i(s_k)=0$ if $i\notin J_n(s_k).$ With
$c_9>0$ being a ''free'' constant, we obtain
\begin{equation}
P_f^n\left( Q_2>c_9\gamma _n\right) \leq \sum_{k=1}^{n^2}P_f^n\left( \frac
1{n\delta _n\rho _{\min }}\left| \sum_{i\in J_n(s_k)}\widetilde{L}_i(s_k)%
\overline{U}(X_i)\right| >c_9\gamma _n\right) .  \label{ES.9}
\end{equation}
Set for brevity $v_n=\widetilde{L}_i(s_k)\left( \rho _{\min }n\delta
_n\gamma _n\right) ^{-1}\log n.$ Then Chebyshev's inequality and the
independence of the r.v.'s $\overline{U}(X_i),$ $i\in I_k$ imply
\begin{equation}
P_f^n\left( \frac 1{n\delta _n\rho _{\min }}\sum_{i\in J_n(s_k)}\widetilde{L}%
_i(s_k)\overline{U}(X_i)>c_9\gamma _n\right) \leq e^{-c_9\log n}\prod_{i\in
J_n(s_k)}E_f^n\exp \left\{ v_n\overline{U}(X_i)\right\} .  \label{ES.10}
\end{equation}
Since $v_n\leq k_{\max }\rho _{\min }^{-1}\kappa _0^{-(2\beta +1)/\beta
}\gamma _n\leq \varepsilon _0$ (for $n$ large enough), by Proposition \ref{P
BKG1}, we obtain
\begin{equation}
\prod_{i\in J_n(s_k)}E_f^n\exp \left\{ v_n\overline{U}(X_i)\right\} \leq
\exp \left\{ \frac{I_{\max }}2v_n^22\tau n\delta _n\right\} \leq \exp
\left\{ c_{10}\log n\right\} ,  \label{ES.11}
\end{equation}
for some constant $c_{10}$ depending on $I_{\max },$ $\kappa _0,$ $\beta ,$ $%
k_{\max },$ $\rho _{\min },$ $\tau .$ Choosing $c_9$ to be $c_{10}+3$ we
get from (\ref{ES.10}) and (\ref{ES.11})
\begin{equation}
P_f^n\left( \frac 1{n\delta _n\rho _{\min }}\sum_{i\in J_n(s_k)}\widetilde{L}%
_i(s_k)\overline{U}(X_i)>c_9\gamma _n\right) \leq \exp \left\{ -3\log
n\right\} .  \label{ES.12}
\end{equation}
In the same way we establish that
\begin{equation}
P_f^n\left( \frac 1{n\delta _n\rho _{\min }}\sum_{i\in J_n(s_k)}\widetilde{L}%
_i(s_k)\overline{U}(X_i)<-c_9\gamma _n\right) \leq \exp \left\{ -3\log
n\right\} .  \label{ES.14}
\end{equation}
From (\ref{ES.12}), (\ref{ES.14}) and (\ref{ES.9}) we get
\begin{equation}
P_f^n\left( Q_2>c_9\gamma _n\right) \leq 2n^2\exp \left\{ -3\log n\right\}
=2\frac 1n.  \label{ES.15}
\end{equation}
Now (\ref{ES.8d}) and (\ref{ES.15}) give us
\begin{equation}
P_f^n\left( \left\| g_n^{*}-E_f^ng_n^{*}\right\| _\infty \geq c_{11}\gamma
_n\right) \leq 4\frac 1n  \label{ES.16}
\end{equation}
for an appropriate constant $c_{11}.$ Finally (\ref{ES.8}) and (\ref{ES.16})
imply (\ref{ES.5}).

Generally speaking $g_n^{*}$ is not bounded. But it is easy to define
another estimator on its basis  which satisfies this requirement. For this
it is enough to put
\[
g_n^{**}(t)=\max \left\{ \min \left\{ g_n^{*}(t),\Lambda _{\max }\right\}
,\Lambda _{\min }\right\}
\]
where $\Lambda _{\min }$ and $\Lambda _{\max }$ are the ends of the interval
$\Lambda _0.$ The estimator $g_n^{**}$ satisfies (\ref{ES.5}) since  for
any $f\in \Sigma $ we have $\Lambda _{\min }\leq g(t)=b(f(t))\leq \Lambda
_{\max },$ which in turn implies
\[
\left\{ \left\| g_n^{**}-g\right\| _\infty >c_3\gamma _n\right\} =\left\{
\left\| g_n^{*}-g\right\| _\infty >c_3\gamma _n\right\} .
\]

An estimator for $f$ can be defined by setting $f_n^{*}(t)=a(g_n^{**}(t)),$
where $a(\lambda )$ is the inverse of the function $b(\theta )$ (see also
Section \ref{sec:VS}). Since the function $a(\lambda )$ is Lipschitz we
obtain from (\ref{VS.3a}),
\[
\left| f_n^{*}(t)-f(t)\right| \leq I_{\min }^{-1}\left|
g_n^{**}(t)-g(t)\right| .
\]
This implies that (\ref{ES.5}) is also satisfied, with $f_n^{*}$ and $f$
replacing $g_n^{*}$ and $g,$ but with other constants (depending also on $%
I_{\min }$). Now we will define an estimator taking a finite number of
values in $\Sigma .$ Since the set $\Sigma $ is compact in the uniform
metric (it is equicontinuous and bounded and hence compact by
Ascoli's theorem), we can cover it by a finite number of balls of
radius $\gamma _n$ and with centers $f_i\in \Sigma ,$ $i=1,...,M.$ The
estimate $\widehat{f}_n$ can be defined to be the $f_i$ closest to the
estimate $f_n^{*}.$ In case of nonuniqueness take the $f_i$ with lowest
index. The estimator constructed has the properties claimed. {%
\mbox{\
\rule{.1in}{.1in}}}

In particular, if we take ${\cal E}$ to be the Gaussian shift experiment,
then  from Lemma \ref{L ES1} we get the following.

\begin{lemma}
\label{L ES2}Let $\beta \in (\frac 12,1).$ In the experiment ${\cal Y}^n$
there is an estimator $\widehat{f}_n:X^n\rightarrow \Sigma $ taking a finite
number of values and fulfilling
\[
\sup_{f\in \Sigma }P_f^n\left( \left\| \widehat{f}_n-f\right\| _\infty
>c_1\gamma _n\right) \leq c_2\frac 1{\sqrt{n}},
\]
where $c_1$ and $c_2$ are constants depending only on $I_{\max },$ $I_{\min
},$ $\kappa _0,$ $L,$ $\beta .$
\end{lemma}


\subsection{Proof of global equivalence}

In this section we prove Theorem \ref{T G2d}.

Let ${\cal E}^n$ and ${\cal Y}^n$ be the experiments defined by (\ref{F.5}),
(\ref{F.6}) and (\ref{G.9}), (\ref{G.10}). Let $f_0\in \Sigma .$ Denote by $%
J^{\prime }$ and $J^{\prime \prime }$ the sets of odd and even numbers,
respectively, in $J=\left\{ 1,...,n\right\} .$ Set
\[
X^{\prime ,n}=\prod_{i\in J^{\prime }}X_{(i)},\quad X^{\prime \prime
,n}=\prod_{i\in J^{\prime \prime }}X_{(i)},\quad R^{\prime \prime
,n}=\prod_{i\in J^{\prime \prime }}R_i,\quad {\bf S}^n=\prod_{i=1}^n{\bf S}%
_i,
\]
where $X_{(i)}=X,$ $R_i=R,$ ${\bf S}_i=X,$ if $i$ is odd, and ${\bf S}_i=R,$
if $i$ is even, $i\in J.$ Consider the following product (local) experiments
corresponding to observations at points $t_i$ with even indexes $i\in J:$
\[
{\cal E}_{f_0}^{\prime \prime ,n}=\bigotimes_{i\in J^{\prime \prime }}{\cal E%
}_{f_0,t_i},\quad {\cal Y}_{f_0}^{\prime \prime ,n}=\bigotimes_{i\in
J^{\prime \prime }}{\cal Y}_{f_0,t_i},
\]
where
\begin{eqnarray*}
{\cal E}_{f_0,t_i} &=&\left( X,{\cal B}(X),\left\{ P_{f(t_i)}:f\in \Sigma
_{f_0}(\gamma _n)\right\} \right) , \\
{\cal Y}_{f_0,t_i} &=&\left( R,{\cal B}(R),\left\{ G_{f(t_i)}:f\in \Sigma
_{f_0}(\gamma _n)\right\} \right) .
\end{eqnarray*}
[cf. (\ref{F.4})-(\ref{F.6}) for the definition of $P_{f(t_i)}$ and (\ref
{G.8})-(\ref{G.9}) for the definition of $G_{f(t_i)}$]. Along with this
introduce global experiments
\[
{\cal E}^{\prime ,n}=\bigotimes_{i\in J^{\prime }}{\cal E}_{t_i},\quad {\cal %
F}^n=\bigotimes_{i=1}^n{\cal F}_i^n,
\]
where ${\cal F}_i^n={\cal E}_{t_i},$ if $i$ is odd, and ${\cal F}_i^n={\cal Y%
}_{t_i},$ if $i$ is even, $i\in J,$ and where
\begin{eqnarray*}
{\cal E}_{t_i} &=&\left( X,{\cal B}(X),\left\{ P_{f(t_i)}:f\in \Sigma
\right\} \right) , \\
{\cal Y}_{t_i} &=&\left( R,{\cal B}(R),\left\{ G_{f(t_i)}:f\in \Sigma
\right\} \right) .
\end{eqnarray*}

We will show that the global experiments ${\cal E}^n$ and ${\cal F}^n$ are
asymptotically equivalent. Toward this end we note that by Theorem \ref{T
VS2d} the experiments ${\cal E}_{f_0}^{\prime \prime ,n}$ and ${\cal Y}%
_{f_0}^{\prime \prime ,n}$ are asymptotically equivalent uniformly in $%
f_0\in \Sigma .$ Theorem \ref{T VS2d} in particular implies that the
one-sided deficiency $\delta \left( {\cal E}_{f_0}^{\prime \prime ,n},{\cal Y%
}_{f_0}^{\prime \prime ,n}\right) $ satisfies for any $f_0\in \Sigma $
\[
\delta \left( {\cal E}_{f_0}^{\prime \prime ,n},{\cal Y}_{f_0}^{\prime
\prime ,n}\right) \leq \varepsilon _n\equiv \left( c_1n^{-\frac{2\beta -1}{%
2\beta +1}}(\log n)^{\frac{14\beta +5}{2\beta +1}}\right) ^{1/2},
\]
where $c_1$ is a constant depending only on $I_{\min },$ $I_{\max },$ $%
\kappa _0,$ $L,$ $\beta .$ Let $\left\| \cdot \right\| $ denote the total
variation norm for measures and let $P_f^{\prime \prime ,n},$ $G_f^{\prime
\prime ,n}$ be the product measures corresponding to the experiments ${\cal E%
}_{f_0}^{\prime \prime ,n}$ and ${\cal Y}_{f_0}^{\prime \prime ,n}:$%
\begin{equation}
P_f^{\prime \prime ,n}=\bigotimes_{i\in J^{\prime \prime }}P_{f(t_i)},\quad
G_f^{\prime \prime ,n}=\bigotimes_{i\in J^{\prime \prime }}G_{f(t_i)}.
\label{PG.2}
\end{equation}
By lemma 9.2 in \cite{Nuss}, for any $f_0\in \Sigma $ there is a Markov
kernel $K_{f_0}^n:(X^{\prime \prime ,n},{\cal B}(R^{\prime \prime
,n}))\rightarrow [0,1]$ such that
\begin{equation}
\sup_{f_0\in \Sigma }\sup_{f\in \Sigma _{f_0}(\gamma _n)}\left\|
K_{f_0}^n\cdot P_f^{\prime \prime ,n}-G_f^{\prime \prime ,n}\right\| \leq
\varepsilon _n.  \label{PG.1}
\end{equation}

Let us establish that there is a Markov kernel $M^n:(X^n,{\cal B}({\bf S}%
^n))\rightarrow [0,1]$ such that
\begin{equation}
\sup_{f\in \Sigma }\left\| M^n\cdot P_f^n-F_f^n\right\| \leq c_2\varepsilon
_n,  \label{PG.2a}
\end{equation}
for some constant $c_2.$ First note that any vector $x\in X^n$ can be
represented as $(x^{\prime };x^{\prime \prime })$ where $x^{\prime }$ and $%
x^{\prime \prime }$ are the corresponding vectors in $X^{\prime ,n}$ and $%
X^{\prime \prime ,n}.$ The same applies for $s\in {\bf S}^n:$ $s=(x^{\prime
};y^{\prime \prime }),$ where $x^{\prime }\in X^{\prime ,n}$ and $y^{\prime
\prime }\in R^{\prime \prime ,n}.$ For any $x=(x^{\prime };x^{\prime \prime
})\in X^n$ and $B\in {\cal B}({\bf S}^n)$ set
\[
M^n(x,B)=\int_{R^{\prime \prime ,n}}{\bf 1}_B\left( (x^{\prime };y^{\prime
\prime })\right) K_{\widehat{f}_n(x^{\prime })}^n(x^{\prime \prime
},dy^{\prime \prime }),
\]
where $\widehat{f}_n(x^{\prime })$ is the preliminary estimator of lemma \ref
{L ES1} in the experiment ${\cal E}^{\prime ,n}.$ It is easy to see that
\begin{eqnarray}
\left( M^n\cdot P_f^n\right) (B) &=&\int_{X^{\prime ,n}}\int_{X^{\prime
\prime ,n}}M^n\left( (x^{\prime };x^{\prime \prime }),B\right) P_f^{\prime
\prime ,n}(dx^{\prime \prime })P_f^{\prime ,n}(dx^{\prime })  \nonumber \\
&=&\int_{X^{\prime ,n}}\int_{R^{\prime \prime ,n}}{\bf 1}_B\left( (x^{\prime
};y^{\prime \prime })\right) \left( K_{\widehat{f}_n(x^{\prime })}^n\cdot
P_f^{\prime \prime ,n}\right) (dy^{\prime \prime })P_f^{\prime
,n}(dx^{\prime })  \label{PG.3}
\end{eqnarray}
and
\begin{equation}
F_f^n(B)=\int_{X^{\prime ,n}}\int_{R^{\prime \prime ,n}}{\bf 1}_B\left(
(x^{\prime };y^{\prime \prime })\right) G_f^{\prime \prime ,n}(dy^{\prime
\prime })P_f^{\prime ,n}(dx^{\prime }),  \label{PG.4}
\end{equation}
where $P_f^{\prime ,n}$ is the measure in the experiment ${\cal E}^{\prime
,n}$ defined analogously  to  $P_f^{\prime \prime ,n}$ in (\ref{PG.2}),
but with $J^{\prime }$ replacing $J^{\prime \prime }.$ By Lemma \ref{L ES1}
there are two constants $c_3$ and $c_4$ depending only on $I_{\max },$ $%
I_{\min },$ $\kappa _0,$ $L,$ $\beta $ such that
\begin{equation}
P_f^{\prime ,n}(A_f^c)\leq c_4\varepsilon _n,  \label{PG.5}
\end{equation}
where $A_f=\left\{ x^{\prime }\in X^{\prime ,n}:\left\| \widehat{f}%
_n(x^{\prime })-f\right\| _\infty \leq c_3\gamma _n\right\} .$ Then (\ref
{PG.3}) and (\ref{PG.4}) imply
\begin{eqnarray*}
\left| \left( M^n\cdot P_f^n\right) (B)-F_f^n(B)\right| &\leq &2P_f^{\prime
,n}(A_f^c) \\
&&+\int_{A_f}\sup_{f_0\in \Sigma }\sup_{f\in \Sigma _{f_0}(\gamma
_n)}\left\| K_{f_0}^n\cdot P_f^{\prime \prime ,n}-G_f^{\prime \prime
,n}\right\| P_f^{\prime ,n}(dx^{\prime }).
\end{eqnarray*}
Using (\ref{PG.1}) and (\ref{PG.5}) we obtain (\ref{PG.2a}). This implies
that the one-sided deficiency $\delta \left( {\cal E}^n,{\cal F}^n\right) $
is less that $c_2\varepsilon _n.$ The bound for $\delta \left( {\cal F}^n,%
{\cal E}^n\right) $ can be obtained in the same way, using Lemma \ref{L ES2}%
. This proves that the Le Cam distance between ${\cal E}^n$ and ${\cal F}^n$
is less that $c_2\varepsilon _n.$ In the same way we can show that ${\cal F}%
^n$ and ${\cal Y}^n$ are asymptotically equivalent. As a result we obtain
asymptotic equivalence of the experiments ${\cal E}^n$ and ${\cal Y}^n.$ As
to the rate of convergence, it is straightforward from the above inequality (%
\ref{PG.2a}) and an analogous one for the pair ${\cal F}^n$ and ${\cal Y}^n.$
Theorem \ref{T G2d} is proved.

{\bf Acknowledgement.} The authors wish to thank Iain Johnstone and Enno
Mammen for pointing out important theoretical connections.

%
%

{\vskip0.5cm}\noindent
{\sc Institute of Mathematics \hfill Weierstrass Institute \newline
Academy of Sciences \hfill Mohrenstr. 39 \newline
Academiei str. 5 \hfill D-10117 Berlin, Germany \newline
Chi\c{s}in\u{a}u 277028, Moldova \hfill e-mail: nussbaum@wias-berlin.de
\newline
e-mail: 16grama@mathem.moldova.su }

\end{document}